\documentclass[11pt]{article}
\usepackage{graphicx}
\usepackage{float}
\usepackage{amsmath}
\usepackage{amsfonts}
\usepackage{comment}
\usepackage{paralist}
\usepackage{amsfonts}
\usepackage{enumerate}
\usepackage{amssymb}
\usepackage{graphicx}
\usepackage{amsthm}
\usepackage{a4wide}
\usepackage{array}
\usepackage{listings}
\usepackage{color}
\usepackage{tikz}
\usetikzlibrary{shapes,arrows,positioning,decorations.pathreplacing}
\usepackage{makecell}
\usepackage{subcaption}
\usepackage{caption}
\usepackage[utf8]{inputenc}
\newcommand {\mat}[1]{\left[\begin{array}{#1}}
\newcommand {\rix}          {\end{array}\right]}
\newcommand {\eq}       [1] {\begin{equation}\label{#1}}
\newcommand {\en}           {\end{equation}}

\newtheorem{theorem}          {Theorem}[section]

\newtheorem{algorithm}         {Algorithm}
\newcommand{\bal}{\begin{center}\begin{tabular}{ll}}
\newcommand{\eal}{\end{tabular}\end{center}}

\catcode`@=11
\font\tenex=cmex10 
\newdimen\p@renwd
\setbox0=\hbox{\tenex B} \p@renwd=\wd0 
\def\bmat#1{\begingroup \m@th
  \setbox\z@\vbox{\def\cr{\crcr\noalign{\kern2\p@\global\let\cr\endline}}%
    \ialign{$##$\hfil\kern2\p@\kern\p@renwd&\thinspace\hfil$##$\hfil
      &&\quad\hfil$##$\hfil\crcr
      \omit\strut\hfil\crcr\noalign{\kern-\baselineskip}%
      #1\crcr\omit\strut\cr}}%
  \setbox\tw@\vbox{\unvcopy\z@\global\setbox\@ne\lastbox}%
  \setbox\tw@\hbox{\unhbox\@ne\unskip\global\setbox\@ne\lastbox}%
  \setbox\tw@\hbox{$\kern\wd\@ne\kern-\p@renwd\left[\kern-\wd\@ne
    \global\setbox\@ne\vbox{\box\@ne\kern2\p@}%
    \vcenter{\kern-\ht\@ne\unvbox\z@\kern-\baselineskip}\,\right]$}%
  \null\;\vbox{\kern\ht\@ne\box\tw@}\endgroup}
\catcode`@=12
\newtheorem{remark} [theorem] {Remark}
\graphicspath{{images/}}
\begin{document}

\title{Error Analysis of a Model Order Reduction Framework for Financial Risk Analysis}
\date{}
\author{Andreas Binder \footnote{MathConsult GmbH, Altenbergerstraße 69, A - 4040 Linz, Austria}
        \and Onkar Jadhav \footnote{Institut f\"ur Mathematik, MA 4-5, TU Berlin, Str. des 17. Juni 136, D-10623 Berlin, Germany}
        \and Volker Mehrmann \footnote{Institut f\"ur Mathematik, MA 4-5, TU Berlin, Str. des 17. Juni 136, D-10623 Berlin, Germany}}
\maketitle

\section*{Abstract}
A parametric model order reduction (MOR) approach for simulating the high dimensional models arising in financial risk analysis is proposed on the basis of  the proper orthogonal decomposition (POD) approach to generate small model approximations for the high dimensional parametric convection-diffusion reaction partial differential equations (PDE). The proposed technique uses an adaptive greedy sampling approach based on surrogate modeling to efficiently locate the most relevant training parameters, thus generating the optimal reduced basis. The best suitable reduced model is procured such that the total error is less than a user-defined tolerance. The three major errors considered are the discretization error associated with the full model obtained by discretizing the PDE, the model order reduction error, and the parameter sampling error. The developed technique is analyzed, implemented, and tested on industrial data of a puttable steepener under the two-factor Hull-White model. The results illustrate that the reduced model provides a significant speedup with excellent accuracy over a full model approach, demonstrating its potential applications in the historical or Monte Carlo value at risk calculations.

\textbf{Keywords:} Financial risk analysis, short-rate models, convection-diffusion-reaction equation, finite element method, parametric model order reduction, proper orthogonal decomposition, adaptive greedy sampling, error analysis, packaged retail investment and insurance-based products.\\

\textbf{MSC(2010):} 35L10, 65M60, 91G30, 91G60, 91G80

\section{Introduction}
Investors use risk analysis to determine whether to undertake a particular venture, the plausible return and how to mitigate potential losses of an activity. The risk analysis of financial instruments often requires the valuation of such instruments under a wide range of future market scenarios. A suitable valuation function takes these market scenarios as input parameters and determines the fair value of financial instruments. Examples of such risk analysis tasks are the calculation of Value-at-Risk (VaR) or the expected shortfall to estimate worst-case scenarios of financial holdings.
Package retail investments and insurance products (PRIIPs) are the packaged instruments that are offered to retail investors. In order to make PRIIPs from different manufacturers more comparable concerning their risks and returns, the European regulation (EU) 1286/2014 requires manufacturers of PRIIPs to supply key information documents (KIDs) to possible retail investors that are easy to read and to understand \cite{EUregKID}. The commission delegated regulation (EU) 2017/653 formulates the details of how the risk and the possible returns of a PRIIP have to be calculated \cite{EUreg}. There are four different categories of PRIIPs. In this paper, we concentrate on category 3 instruments, for which at least $10\,000$ market data scenarios have to be generated, and the PRIIP has to be valuated under these scenarios.

In our approach, the financial instruments are evaluated via the dynamics of short-rate models \cite{Brigo06}, based on convection-diffusion-reaction partial differential equations (PDE). The choice of the short-rate model depends on the underlying financial instrument. Some of the prominent financial models are the one-factor Hull-White model \cite{Hull94}, the shifted Black-Karasinski model \cite{black91}, and the two-factor Hull-White model \cite{Hull01}. These models are calibrated based on market data like \emph{yield curves} that generate a high dimensional parameter space \cite{Brigo06}. To perform the risk analysis, the financial model needs to be solved for such a high dimensional parameter space, and this  requires efficient algorithms. For instruments with high complexity and long-time horizons, computing times of minutes for a single valuation are not unusual.

Thus, in a previous paper \cite{binder2020}, we have established a parametric model order reduction approach based on a variant of the proper orthogonal decomposition approach, which significantly reduces the overall computation time \cite{Berkooz93,Chatterjee00}.
This nonlinear MOR approach is computationally feasible \cite{Liang02} as it determines low dimensional linear (or affine) subspaces \cite{Rathinam03,Vidal69} via a truncated singular value decomposition of a snapshot matrix \cite{Sirovich87}, that is computed by simulating the full model obtained by discretizing the PDE for a small number of pre-selected training parameter values. The question of how to select these parameters is often the most challenging part of the process.
Our previous work \cite{binder2020} has established greedy algorithms to determine the training parameters more efficiently. The adaptive greedy approach searches the entire parameter space efficiently using a surrogate model and determines the best suitable training parameters.
The principal component regression technique has been used to construct the surrogate model.

This paper provides a detailed error analysis of the model order reduction framework. We have considered the hierarchical framework to perform the error analysis that starts by considering the solution of the partial differential equation as the true solution. We then discretize the underlying PDE to generate the full model and project it on a reduced basis to obtain the reduced model. Each stage described here is associated with some numerical error. Three major numerical errors considered are the discretization error associated with the full model, the model order reduction error, and the sampling error associated with the sampling algorithms. The total error associated with the model order reduction framework is the sum of these three significant errors. Our approach aims to minimize the total error such that it is less than the user-defined tolerance generating the suitable reduced model. We define a discretization error estimator based on estimating the exact solution of the PDE, which is higher-order accurate than the underlying numerical solution \cite{Roy10}. The method is known as Richardson extrapolation \cite{Richardson1911}, which uses a sequence of solutions to estimate the discretization error and can be applied to almost any type of discretization approach as well as to both local and global quantities \cite{Roy09}.
The discretization error is mainly governed by the choice of the grid size, so we present an approach to select an optimal grid  based on the developed error estimator.

The model order reduction error is composed of the projection error associated with the proper orthogonal decomposition approach and how well the reduced model approximates the full model \cite{Kunisch01}. Both the projection error and the reduced model error depend on the dimension of the  reduced basis. We present an algorithm to select an optimal reduced dimension that minimizes the projection error and the reduced model error. To obtain a reduced basis, we have to compute a singular value decomposition of the snapshot matrix, which is a computationally costly task. Recent research has shown that the randomized singular value decomposition is a good alternative and it is noticeably faster than the basic SVD for some applications \cite{Mahoney11,oh2015,Song17}. Thus, we use the randomized truncated singular value decomposition, developed in \cite{Halko11}, in the adaptive greedy sampling approach to determine the reduced basis. The sampling algorithm used to obtain the training parameters to construct a reduced basis introduces a sampling error. In the adaptive greedy sampling, the projection error due to the principal component analysis (PCA) and the prediction error of the surrogate model  form the total sampling error. We analyze the projection error due to the PCA in a similar sense to the projection error of the POD, while the mean square error of prediction is used to assess the quality of the surrogate model.

We also perform the sensitivity analysis for these errors to know the relative contribution of each error to the total error \cite{Saltelli10}. It provides useful insight into which  error contributes most to the variability of the model output \cite{Frey02}. The sensitivity analysis assists in improving the model output by selecting the proper grid size, reduced dimension, and training parameters. It also helps to allocate computational resources efficiently as more computational effort could be spent on the most sensitive parameters.

In a previous paper, we have considered a numerical example of a floater with caps and floors under the one-factor Hull-White model \cite{binder2020}. However, there are more complicated instruments, e.g., steepeners whose coupons depend on the difference between the two interest rates \cite{Binder013}. For such instruments, we need to consider a multi-factor model. We solve a numerical example of the puttable steepener with caps and floors using the two-factor Hull-White model in this paper.

\section{Model Order Reduction}
\label{MOR}
In an affine one-factor model, the \emph{forward rates} $f(t,T_1), f(t,T_2)$ at two different tenors $T_1, T_2$ are perfectly correlated, i.e., their correlation coefficient is one, which is unrealistic \cite{Bohner11}. This perfect correlation only allows parallel shifts of the yield curves. However, we often noticed the yield curve steeping where short-term rates are low and the long-term rates are high \cite{Brigo06}. Also, for instruments like steepeners \cite{Albrecher13}, where coupons depend on the difference between two rates, we need to consider a multi-factor model like the two-factor Hull-White model. The two-factor model overcomes these drawbacks  by adding a stochastic disturbance to the drift term.
In the two-factor Hull-White model, the short rate $r$ is assumed to satisfy the following stochastic differential equation \cite{Hull094}
\begin{equation}
\begin{aligned}
    dr(t) &= (\theta(t) + u(t) - \alpha(t)r(t))dt + \sigma_1(t)dW_1(t),\\
    du(t) &= -b(t)u(t) + \sigma_2(t) dW_2(t).
\end{aligned}
\label{1}
\end{equation}
where $\alpha, b, \sigma_1, \sigma_2 > 0$ and $W_1, W_2$ are two Brownian motions under the risk neutral measure such that
\begin{equation*}
    dW_1(t) dW_2(t) = \gamma dt,
\end{equation*}
where $-1 \leq \gamma \leq 1$ is the correlation coefficient.
To evaluate the underlying financial instrument, we use the partial differential equation for the two-factor Hull-White model and then solve it using appropriate boundary conditions and known terminal conditions.
Consider a financial instrument $V(t,r(t),u(t))$ contingent on the stochastic interest rate movement $r(t)$. The two-factor Hull-White PDE is then given as \cite{Hull094}
\begin{equation}
    \frac{\partial V}{\partial t} + (\theta(t) + u - \alpha(t) r) \frac{\partial V}{\partial r} - b(t)u \frac{\partial V}{\partial u} + \frac{1}{2} \sigma_1^2(t) \frac{\partial^2 V}{\partial r^2} + \frac{1}{2} \sigma_2^2(t) \frac{\partial^2 V}{\partial u^2} + \gamma \sigma_1(t) \sigma_2(t) \frac{\partial V}{\partial r \partial u} - rV = 0
    \label{2}
\end{equation}
and the deterministic function $\theta(t)$ is chosen to fit the simulated yield curves. However, the regulation does not provide a methodology for obtaining the parameters $\alpha(t),b(t),\sigma_1(t),\sigma_2(t)$. Thus, considering them as constants should lead to more robust results than time-dependent parameters \cite{Binder013}. We consider the resulting \textit{robust two-factor Hull-White model} with parameters $\alpha,b,\sigma_1,\sigma_2 > 0$ as positive constants and only the parameter $\theta(t)$ as time-dependent. Our results can, however, be extended to the more general case.
According to the PRIIPs regulations, we have to perform at least $s= 10\, 000$ yield curve simulations. We construct a simulated yield curve matrix
\begin{equation*}
    Y =
\begin{bmatrix}
y_{11} & \cdots & y_{1m}\\
\vdots & \vdots & \vdots\\
y_{s1} & \cdots & y_{sm}\\
\end{bmatrix} \in \mathbb{R}^{s\times m},
\end{equation*}
which is then used to calibrate the parameter $\theta(t)$.
The calibration based on $s=10\,000$ different simulated yield curves generates $s$ different piecewise constant parameters $\theta_\ell(t)$, which change their values $\theta_{\ell,i}$ only at the $m$ tenor points. We incorporate these in a  matrix
\begin{equation*}
    \Theta=
\begin{bmatrix}
\theta_{11} & \cdots & \theta_{1m}\\
\vdots & \vdots & \vdots\\
\theta_{s1} & \cdots & \theta_{sm}\\
\end{bmatrix}.
\end{equation*}
A detailed procedure for the yield curve simulation and the calibration of $\theta(t)$ based on simulated yield curves is described in \cite{binder2020,MathConsult09}, respectively.
We have applied a finite element method to solve the PDE numerically, for more details see \cite{Binder013}. This discretization is a parametric high dimensional model of the  form
\begin{equation}
    A(\rho_\ell(t))V^{n-1} = B(\rho_\ell(t))V^{n},
    \label{3}
\end{equation}
with given terminal vector $V^T$, and matrices $A(\rho_\ell) \in \mathbb{R}^{M\times M}$, and $B(\rho_\ell) \in \mathbb{R}^{M\times M}$. Here $\rho = \{ \alpha,b,\sigma_1,\sigma_2,\gamma, \theta(t)\}$ is the group of model parameters. We call this the \emph{full model} (FM) for the model reduction procedure. We solve (\ref{3}) by propagating backward in time. Here again $\ell = 1,\dots,s=10\, 000$,  $m$ is the total number of tenor points, and we need to solve this system at each time step $n$ with an appropriate boundary conditions and a known terminal value for the underlying instrument.
Altogether we have a parameter space $\mathcal P$ of size $10\, 000 \times m$ to which we now apply model reduction.

To perform the parametric model reduction for system (\ref{3}), we employ Galerkin projection onto a low dimensional subspace via
\begin{equation}
    \Bar{V}^n = QV^n_d,
    \label{4}
\end{equation}
where the columns of $Q \in \mathbb{R}^{M\times d}$ represent the reduced basis with $d \ll M$, $V_d^n$ is a vector of reduced coordinates, and $\Bar{V}^n \in \mathbb{R}^{M}$ is the solution in the $n$th time step obtained using the reduced model. For the Galerkin projection we require that the residual of the reduced state
\begin{equation}
    p^n(V_d^n,\rho_\ell) = A(\rho_\ell)QV^{n-1}_d - B(\rho_\ell)QV^n_d,
    \label{5}
\end{equation}
is orthogonal to the reduced basis matrix $Q$, i.e.,
\begin{equation}
    Q^Tp^n(V^n_d,\rho_\ell) = 0,
    \label{6}
\end{equation}
so that by multiplying $p^n(V_d^n,\rho_\ell)$ with  $Q^T$, we get
\begin{equation}
\begin{aligned}
    Q^TA(\rho_\ell)QV^{n-1}_d &= Q^TB(\rho_\ell)QV^n_d,\\
    A_d(\rho_\ell)V^{n-1}_d &= B_d(\rho_\ell)V^n_d,
    \label{7}
\end{aligned}
\end{equation}
where $A_d(\rho_\ell) \in \mathbb{R}^{d\times d}$ and $B_d(\rho_\ell) \in \mathbb{R}^{d \times d}$ are the parameter dependent reduced matrices. 

We obtain the Galerkin projection matrix $Q$ in (\ref{4}) based on a proper orthogonal decomposition (POD) approach, which generates an optimal order orthonormal basis $Q$ in the least square sense that is independent of the parameter space $\mathcal{P}$ and we do this by the method of snapshots.
The snapshots are nothing but the state solutions obtained by simulating the full model for selected parameter groups. We assume that we have a training set of parameter groups $\rho_{1}, \ldots,\rho_{k}\in [\rho_{1},\rho_{s}]$. We compute the solutions of the full model for this training set and combine them in a snapshot matrix $\hat{V} = [V(\rho_1),V(\rho_2),...,V(\rho_k)]$. The POD method solves, see \cite{Golub70},
\begin{equation*}
    \mathrm{POD}(\hat{V}) := \underset{Q}{\mathrm{argmin}} \frac{1}{k} \sum_{j=1}^{k} \|V_j - QQ^TV_j\|^2,
\end{equation*}
for an orthogonal matrix $Q \in \mathbb{R}^{M \times d}$ via a \emph{truncated singular value decomposition (SVD)}
\begin{equation*}
   \hat{V} = \Phi \Sigma \Psi^T=\sum_{i=1}^k \Sigma_i\phi_i\psi_i^T,
\end{equation*}
where $\phi_i$ and $\psi_i$ are the left and right singular vectors of the matrix $\hat{V}$ respectively, and $\Sigma_i$ are the singular values.  The truncated  SVD computes only the first $k$ columns of the matrix $\Phi$. The detailed error analysis of the model order reduction and the randomized truncated singular value decomposition is explained in Subsection \ref{MORError}.
We choose only $d$ out of $k$ POD modes to construct $Q = [\phi_1 \cdots \phi_{d}]$ which minimizes the projection error (\ref{28})
\begin{equation*}
    \epsilon_{\mathrm{POD}} = \frac{1}{k} \sum_{j=1}^{k} \| V_j - \sum_{k=1}^{\ell}(V_j \phi_k) \phi_k\|^2 = \sum_{\ell = d+1}^{k} \Sigma_\ell^2.
\end{equation*}
We summarize the procedure of selecting the dimension $d$ of the reduced basis in Algorithm~\ref{Algo5}.
It is evident that the quality of the reduced model strongly depends on the selection of parameter groups $\rho_1,...,\rho_k$ that are used to compute the snapshots. Hence it is essential to introduce an efficient sampling technique for the  parameter space. We could consider standard sampling techniques, like uniform sampling or random sampling \cite{James013}. However, these techniques may neglect vital regions within the parameter space. As an alternative, a greedy sampling method has been suggested in the framework of model order reduction.
The greedy sampling technique selects the parameters that maximize the error between the reduced model and the full model. Further, the reduced basis is obtained using these selected training parameters. The relative error calculation is expensive as it requires the full model solution, so instead, the residual error associated with the reduced model (\ref{5}) is used as an error estimator. However, it is not reasonable to compute an error estimator for the entire parameter space. This means that one has to select a pre-defined parameter set as a subset of the high dimensional parameter space to train the greedy sampling algorithm. We usually select this pre-defined subset randomly. A random selection may neglect the crucial parameters within the parameter space. Thus, to surmount this problem, we implemented an adaptive greedy sampling approach. The algorithm chooses the most suitable parameters adaptively using an optimized search based on surrogate modeling. This approach evades the cost of computing the error estimator for each parameter within the parameter space and instead uses a surrogate model to locate the training parameter set. We have built the surrogate model using the principal component regression technique. See \cite{binder2020} for detailed sampling algorithms and their implementation. This paper aims to study the numerical errors that arise during different stages of the model order reduction approach.
\section{Error Analysis}
\label{ErrAnalysis}
In most cases, the exact solutions of mathematical models are unknown, and they incorporate some modeling errors. We have no other way to 
compute approximate solutions and to analyze them. However, the approximate solution contains some numerical errors. It is necessary to reduce these numerical errors as much as possible to improve the level of accuracy of the numerical solutions. Our  mathematical model  is a partial differential equation (PDE) with appropriate boundary conditions and a known terminal condition. To solve this system, we have designed an approximate solution by discretizing the PDE and thus committing a discretization error. Using model reduction based on the POD leads to a  projection error, which depends on the excluded singular values while generating the reduced basis. To know how well the reduced model approximates the full model, we have to determine the model order reduction error. In this work, the training parameters are chosen based on either the classical greedy sampling or the adaptive greedy sampling approaches. Thus, we also have to consider a sampling error associated with the sampling algorithms as well.

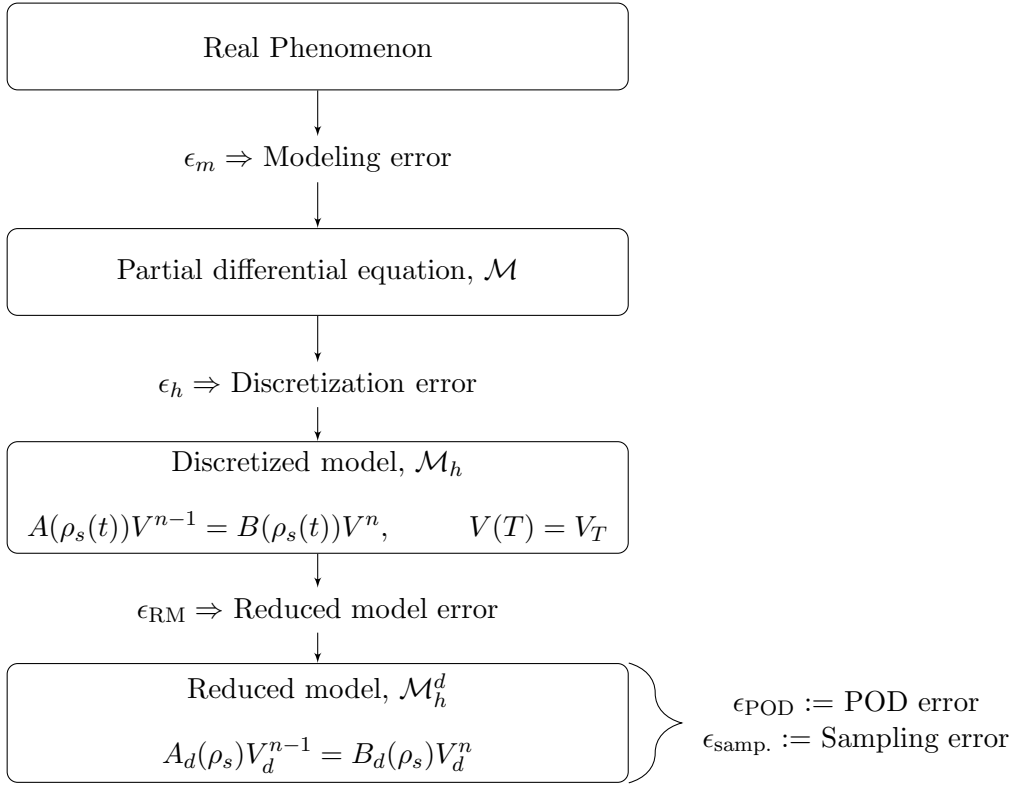
\begin{figure}[htb]
\begin{center}
\usetikzlibrary{shapes,arrows,positioning,decorations.pathreplacing}
\tikzstyle{block} = [rectangle, draw,
    text width=8cm, text centered, rounded corners, minimum height=3em]
\tikzstyle{line} = [draw, -latex']

\begin{tikzpicture}[node distance = 1.5cm, auto]
\node[block](Model){Real Phenomenon};
\node[draw=none,fill=none, below of=Model] (ModlError){$\epsilon_m \Rightarrow$ Modeling error};
\node [block, below of=ModlError] (MainModel) {Partial differential equation, $\mathcal{M}$
};
\node[draw=none,fill=none, below of=MainModel] (DisError){$\epsilon_h \Rightarrow$ Discretization error};
\node[block, below of=DisError](DiscModel){Discretized model, $\mathcal{M}_h$
\begin{equation*}
    A(\rho_s(t))V^{n-1} = B(\rho_s(t))V^{n}, \hspace{1cm} V(T) = V_T
\end{equation*}
};
\node[draw=none,fill=none, below of=DiscModel](MORError){$\epsilon_{\mathrm{RM}} \Rightarrow$ Reduced model error};
\node[block, below of=MORError](RM){Reduced model, $\mathcal{M}_h^d$
\begin{equation*}
    A_d(\rho_s)V^{n-1}_d = B_d(\rho_s)V^n_d
\end{equation*}
};
\path [line] (Model) -- (ModlError);
\path [line] (ModlError) -- (MainModel);
\path [line] (MainModel) -- (DisError);
\path [line] (DisError) -- (DiscModel);
\path [line] (DiscModel) -- (MORError);
\path [line] (MORError) -- (RM);
\draw[decorate,decoration={brace,amplitude=20pt}]
(RM.north east) -- (RM.south east) node [midway,xshift=5.2cm,left,align=center] {$\epsilon_{\mathrm{POD}} :=$ POD error \\$\epsilon_{\mathrm{samp.}} :=$ Sampling error};
\end{tikzpicture}
\end{center}
\caption{A model hierarchy showing errors arising in the analysis of the mathematical model.}
\label{fig:1}
\end{figure}

Figure \ref{fig:1} shows the model hierarchy to obtain the reduced model with associated errors in each stage.
Usually, the formulated mathematical model is not an exact representation of the underlying phenomenon but subject to  many simplifying assumptions.  In this paper, we assume the solution of the partial differential equation as the \emph{true solution}, i.e., the modeling error $\epsilon_m$ will be neglected.

Let $V$ be the exact solution for the mathematical model $\mathcal{M}$.
Let $V_h$ be the solution obtained using the discretized model $\mathcal{M}_h$ with the mesh size $h$. Then the approximate solution $V_h$ encompasses the discretization error
\begin{equation}
    \epsilon_h = \frac{\| V - V_h \|}{\|V\|}.
    \label{8}
\end{equation}
Furthermore, the reduced model $\mathcal{M}^d_h$ of dimension $d$ is obtained. Let the $\Bar{V}$ be the result obtained from the reduced model. The reduced model error can then be given as
\begin{equation}
    \epsilon_\mathrm{RM} =  \frac{\| V_h - \Bar{V} \|}{\|V_h\|}.
    \label{9}
\end{equation}
The result obtained using the reduced model encloses the total numerical error, and based on that, we can write the final output equation with numerical errors
\begin{equation}
    \|V - \Bar{V}\| = \|V - V_h\| + \| V_h - \Bar{V}\|.
    \label{10}
\end{equation}
The reduced model is obtained based on the proper orthogonal decomposition approach, which induces a projection error $\epsilon_{\mathrm{POD}}$. The total model order reduction error $\epsilon_d$ is then the sum of a projection error $\epsilon_{\mathrm{POD}}$, and a reduced model error $\epsilon_{\mathrm{RM}}$
\begin{equation*}
    \epsilon_d \approx \epsilon_{\mathrm{RM}} + \epsilon_{\mathrm{POD}}.
\end{equation*}
We use a sampling approach to select the most relevant parameters. This produces a reasonably accurate reduced basis, but the sampling algorithms contribute to an additional sampling error. Thus, the total simulation error $\epsilon_T$ can be estimated by the discretization error $\epsilon_h$, a sampling error $\epsilon_{\mathrm{samp}}$ and the model order reduction error $\epsilon_d$ as
\begin{equation}
    \epsilon_T \approx \epsilon_m + \epsilon_h + \epsilon_{\mathrm{RM}} + \epsilon_{\mathrm{POD}} + \epsilon_{\mathrm{samp}}.
    \label{11}
\end{equation}
We aim to generate the reduced model and to solve the financial instrument for the entire parameter space such that the as a cost function the total error is minimized. This total error depends on the grid size $h$, the reduced dimension $d$, and the quality of training parameters $c$ used to obtain the reduced basis.
We optimize these variables 
such that the total error $\epsilon_T$ is less than the user-defined tolerance $e_{tol}$, i.e.,
\begin{equation}
    \underset{\{h,d,c\}}{\mathrm{min}} \epsilon_h + \epsilon_d + \epsilon_{\mathrm{samp}},
    \label{12}
\end{equation}
\begin{equation*}
    \epsilon_T \leq e_{tol}.
\end{equation*}
Another important issue is to know the contribution and the effect of each defined error on the total error (\ref{11}).
Sensitivity analysis enables identifying the errors that have the most significant influence on the model output. 
The sensitivity approach presented in Subsection \ref{Sensitivity} ranks the numerical errors according to their contribution to the total error and helps to efficiently allocate the computational resources. The extra computational effort can be spent on the cost parameters most sensitive to the final output at the expense of the computational cost required to obtain the least sensitive parameters.
\subsection{Discretization of a PDE}
\label{DiscPDE}
For the discretization of a PDE, we use a finite dimension grid $\bigg (\prod_{k=1}^{sd} [r_k,u_k] \bigg ) \times [0, T]$ in $sd$ spatial dimensions and one temporal dimension. The simulation time $[0, T]$ is divided into $N$ equal intervals, $[t^n, t^{n-1}]$ where $n=1,\cdots,N$. We have considered an implicit time-stepping scheme. To obtain an optimal solution of the discrete model, one must determine the best suitable time step and the grid size. We implement an adaptive time-stepping scheme to achieve the optimal time step, which solves the model at the desired accuracy but with low computational cost while the grid size is selected such that the discretization error is less than a specified user-defined tolerance.
We have a PDE of the form
\begin{equation*}
    \frac{\partial V}{\partial t} = L(r,u) V,
\end{equation*}
where $L(r,u)$ is a linear differential operator containing all spatial derivatives.
The idea is to obtain the optimal time step $\Delta t_{opt}$ such that it minimizes the relative error between the exact solution $V$ and the solution obtained using the time step $\Delta t_{opt}$
\begin{equation*}
    \epsilon_t = \|V - V_h(\Delta t_{opt}) \|.
\end{equation*}
Assuming the numerical method used is of order $p$ in time, the exact value of the PDE $V$, and the approximated value $V_h(\Delta t)$ with the time step $\Delta t$, the error for the approximation is given as \cite{Soderlind02}
\begin{equation*}
    V = V_{h}(\Delta t) - C \Delta t^{p+1} - \mathcal{O}(\Delta t^{p+2}) .
\end{equation*}
By neglecting the higher order terms we can approximately determine the value of the constant $C$ as
\begin{equation}
    C = \frac{V_h(\Delta t) - V}{\Delta t^{p+1}}.
    \label{13}
\end{equation}
Let $V_h(K\Delta t)$ be the solution obtained using the time step $K\Delta t$
\begin{equation}
    V = V_{h}(K\Delta t) - C K^{p+1}\Delta t^{p+1} - \mathcal{O}(\Delta t^{p+2}).
    \label{14}
\end{equation}
Considering a time discretization scheme of order $p=1$ and substituting (\ref{14}) into (\ref{13}), we get
\begin{equation}
    C = \frac{V_{h}(K\Delta t) - V_{h}(\Delta t)}{\Delta t^2 (K^2 -1)}.
    \label{15}
\end{equation}
We calculate $\Delta t_{opt}$ such that the relative error $\epsilon_t$ is less than the user-defined tolerance
\begin{equation*}
    \epsilon_t = \| V - V_h(\Delta t_{opt}) \| \leq e_{tol}^t,
\end{equation*}
where
\begin{equation*}
    V_{h}(\Delta t_{opt}) = V + C \Delta t_{opt}^2 + \mathcal{O}(\Delta t_{opt}^3).
\end{equation*}
Thus, we get the relative error estimate $\hat{\epsilon}_t$ by substituting $C$ as
\begin{equation*}
    \|V - V_h(\Delta t_{opt}) \| \approx \bigg (\frac{\Delta t^2_{opt}}{\Delta t^2} \bigg ) \frac{\| V_h(\Delta t) - V_h(K\Delta t) \|}{K^2 - 1} = \hat{\epsilon}_t.
\end{equation*}
To determine the optimal step size we set
\begin{equation*}
    \hat{\epsilon}_t = e_{tol}^t,
\end{equation*}
which yields
\begin{equation}
    \Delta t_{opt}^2 \approx \frac{e_{tol}^t \Delta t^2 (K^2 - 1)}{\| V_h(\Delta t) - V_h(K\Delta t) \|}.
    \label{16}
\end{equation}
Since the roundoff or iterative convergence errors that arise while solving the system of linear equations are negligible compared to the discretization error, we can neglect them \cite{Liang11}.
The adaptive time-stepping strategy aims to find the time step that will achieve the desired accuracy at a relatively low cost \cite{Belfort07,Einkemmer18}. Thus, we have to determine the largest possible time step that satisfies (\ref{16}).
An automatic time step control can be executed using the optimal time step $\Delta t_{opt}$, as shown in Algorithm \ref{Algo1}. The step size controller is designed on the assumption that the largest possible step size should be selected. If $\Delta t_{opt} > \Delta t$, we can say that we have a good estimation of the solution that satisfies the user-defined tolerance. Thus, the algorithm runs until the optimal step size $\Delta t_{opt}$ becomes larger than the step size $\Delta t$.
\begin{algorithm}[\textbf{Automatic time step control}]
~\\
\textbf{Input:} Full model $\mathcal{M}_h$, maximum time step $\Delta t_{max}$, $e_{tol}^t$, $K$. \\
\textbf{Output:} $\Delta t_{opt}$.
\begin{itemize}
    \item [ ] Solve the full model with the time step $\Delta t_{max}$ and store the result in $V_h(\Delta t_{max})$.
    \item [ ] WHILE $\Delta t_{opt} < \Delta t$
        \item [ ] \quad Compute a smaller step size $\Delta t = \Delta t_{max}/K$.
        \item [ ] \quad Solve the full model with the time step $\Delta t$ and store the result in $V_h(\Delta t)$.
        \item [ ] \quad Compute an optimal time step $\Delta t_{opt}= \sqrt{\frac{e_{tol}^t \Delta t^2 (K^2 - 1)}{\| V_h(\Delta t) - V_h(\Delta t_{max}) \|}}$.
        \item [ ] \quad Increment $K = \mathrm{round}(\Delta t/\Delta t_{opt}) + 1$.
    \item [ ] END.
    \item [ ] Set $\Delta t_{opt} = \Delta t$.
\end{itemize}
\label{Algo1}
\end{algorithm}
It is important to note that we have 
to select the time step such that all key time points are achievable, i.e., coupon dates, put dates, or valuation dates.

We also need a discretization error estimator based on estimating the exact solution for the differential equation, see 
\cite{Roy10}. We use the Richardson extrapolation approach that estimates the  exact solution using the formal rate of convergence of a discretization method and two or three different solutions on systematically refined meshes \cite{Roy09}. One would argue that the computation of multiple simulations in the asymptotic range of the full model increases the computational burden. However, the main advantage of the Richardson approach is that it can be used as a post-processing technique and can be applied to any discretization method (e.g., finite element or finite difference). It estimates the total error, including locally generated errors and those transported from other regions of the domain, see also \cite{Almohammadi13,Kwasniewski13,Phillips14}. 

Consider a discretization error as
\begin{equation}
    \epsilon_h =\bigg| \frac{ V_h - V}{V} \bigg|,
    \label{17}
\end{equation}
where $V$ is the exact solution of the mathematical model and $V_h$ is a solution obtained on a grid spacing $h$ using a discretized model. In the case of a two-dimensional system
\begin{equation*}
    h = \bigg[ \frac{1}{E} \sum_{i=1}^E A^{(e)}_i \bigg]^{1/2},
\end{equation*}
where $E$ is the number of elements, and $A^{(e)}_i$ is the area of $i$th element.
For a smooth enough solution with no discontinuities, we can expand the numerical solution $V_h$ using Taylor  expansion
\begin{equation*}
    V_h = V + \xi_p h^p + \xi_{p+1}h^{p+1} + \xi_{p+2} h^{p+2} + \cdots,
\end{equation*}
which leads to the discretization error
\begin{equation*}
    \epsilon_h = |V_h - V| = \xi_p h^p + \xi_{p+1}h^{p+1} + \xi_{p+2} h^{p+2} + \cdots = \sum_{p=P_f} \xi_p h^p \approx \xi_{p_f} h^{p_f} :=  \Bar{\epsilon}_h,
\end{equation*}
by neglecting higher order terms. Here $p_f$ is the formal order of accuracy determined by the preferred discretization scheme and is the first term in the series. 
To accurately estimate the discretization error $\Bar{\epsilon}_h$, 
Richardson extrapolation demands that all other sources of numerical errors (e.g., round-off error, iterative convergence errors) are negligible compared to the discretization error \cite{Liang11}.

Consider a grid refinement index $g$ as the ratio of coarse to fine mesh size. It typically is recommended  see e.g. \cite{Celik08}, to use the grid refinement ratio between $1.5$ to $2$ as it allows to distinguish the discretization error from other error sources.
\begin{equation*}
    g = \frac{h_{\mathrm{coarse}}}{h_{\mathrm{fine}}} > 1.
\end{equation*}
Let $h_{\mathrm{fine}} = h$, and $h_{\mathrm{coarse}} = gh$, then the discretization error estimator is \cite{Roy10}
\begin{equation}
    \Bar{\epsilon}_h = \frac{|V_h - V|}{V} = \frac{1}{(g^p - 1)} \bigg| \frac{V_h - V_{gh}}{V_h} \bigg|.
    \label{18}
\end{equation}
To access the confidence in the discretization error estimate, it is necessary to define the observed order of accuracy $\hat{p}$ \cite{Roy10}.
The requirement of the reliability of the discretization error estimate is that the solutions are in the range of asymptotic convergence. To demonstrate that the asymptotic range for different solutions has been achieved, we must have at least three different discretized solutions with three different mesh sizes. We further calculate the observed order using these discretized solutions over a range of meshes.
According to \cite{Celik08}, when the observed order of accuracy is equal to the formal order of accuracy, we can have a high degree of confidence in the error estimator.
We can determine the observed order of accuracy as follows.
Consider three different mesh sizes $h_1$, $h_2$ and $h_3$ such that
\begin{equation*}
    g_{12} = \frac{h_2}{h_1} > 1 \hspace{1cm} g_{23} = \frac{h_3}{h_2} >1.
\end{equation*}
Let $h_1 = h$ then $h_2 = g_{12}h$ and $h_3 = g_{23}g_{12}h$.
Using the formal power series
\begin{align}
    V_1 &= V + \xi_p h^p + \xi_{p+1} h^{p+1} + \cdots, \label{19}\\
    V_2 &= V + \xi_p g_{12}^p h^p + \xi_{p+1}g_{12}^{p+1} h^{p+1} + \cdots, \label{20} \\
    V_3 &= V + \xi_p (g_{23} g_{12})^p h^p + \xi_{p+1}(g_{23}g_{12})^{p+1} h^{p+1} + \cdots, \label{21}
\end{align}
and subtracting (\ref{20}) from (\ref{21}) and (\ref{19}) from (\ref{20}) and setting $\hat{p} \approx p$,
we obtain
\begin{align}
    V_3 - V_2 &= \xi_{\hat{p}} h^{\hat{p}}(g_{23}^{\hat{p}} g_{12}^{\hat{p}} - g_{12}^{\hat{p}}) + \mathcal{O}(h^{\hat{p} +1}), \label{22}\\
    V_2 - V_1 &= \xi_{\hat{p}} h^{\hat{p}}(g_{12}^{\hat{p}} -1) + \mathcal{O}(h^{\hat{p} +1}). \label{23}
\end{align}
Dividing (\ref{22}) by (\ref{23}) and neglecting higher order terms we have
\begin{align*}
    \frac{V_3 - V_2}{V_2 - V_1} &= \frac{g_{23}^{\hat{p}} g_{12}^{\hat{p}} - g_{12}^{\hat{p}}}{g_{12}^{\hat{p}} -1},\\
    \frac{V_3 - V_2}{g_{23}^{\hat{p}} - 1} &= g_{12}^{\hat{p}}\bigg ( \frac{V_2 - V_1}{g_{12}^{\hat{p}} -1} \bigg ).
\end{align*}
We can then determine $\hat{p}$ by an iterative process 
\begin{equation}
    \hat{p}_{k+1} = \frac{\mathrm{ln} \Bigg[( {g_{12}}^{\hat{p}_k} - 1) \bigg( \frac{V_3 - V_2}{V_2 - V_1} \bigg) + {g_{12}}^{\hat{p}_k} \Bigg]}{\mathrm{ln}(g_{12}g_{23})},
    \label{24}
\end{equation}
where an initial guess for $\hat{p}$ is the formal order of accuracy. Algorithm \ref{Algo2} shows the procedure to compute the observed order of accuracy. We consider it converged when the change in two subsequent $\hat{p}_{k+1}$ and $\hat{p}_{k}$ is not significant.
\begin{algorithm}[\textbf{Algorithm to estimate the observed order of accuracy}]
~\\
\textbf{Input:} $V_1, V_2, V_3$, $\hat{p}_1 = p_f$. \\
\textbf{Output:} $\hat{p}$.
\begin{itemize}
    \item [ ] Set $\hat{p}_1 = p_f$.
    \item [ ] $\Delta V =  (V_3 - V_2)/(V_2 - V_1)$.
    \item [ ] k = 1.
    \item [ ] WHILE $\Delta \hat{p} < 10^{-3}$
    \item [ ] \quad Compute (\ref{24}): $ \hat{p}_{k+1} = \mathrm{ln} \bigg[( {g_{12}}^{\hat{p}_k} - 1) \Delta V + {g_{12}}^{\hat{p}_k} \bigg] \big/ \mathrm{ln}(g_{12}g_{23})$.
    \item [ ] \quad $\Delta p = \hat{p}_{k+1} - \hat{p}_k$.
    \item [ ] END
    \item [ ] Set $\hat{p} = \hat{p}_k$
\end{itemize}
\label{Algo2}
\end{algorithm}
Substituting $p= \hat{p}$ and $g= g_{12}$ in (\ref{18}) leads to a discretization error estimate \cite{Roy10}
\begin{equation}
    \Bar{\epsilon}_h = \frac{1}{(g_{12}^{\hat{p}} - 1)} \bigg| \frac{V_1 - V_2}{V_1} \bigg|.
    \label{25}
\end{equation}
When the observed order of accuracy $\hat{p}$ matches the formal order of accuracy $p_f$, we have high confidence that this error estimate is accurate. However, it may happen that the observed order of accuracy does not match the formal order of accuracy \cite{Phillips14}, then the error estimator is not  reliable. This may be due to the failure of obtaining the solutions within the asymptotic range or the lack of additional information. Also, the implemented boundary conditions and the terminal/initial conditions change the observed order of accuracy. However, we can fix this problem by converting the error estimator into an epistemic uncertainty estimator \cite{Roy10}. Epistemic uncertainties are different from random or aleatory uncertainties, as epistemic uncertainties arise due to the lack of knowledge.

A grid convergence index (GCI) measures how far the computed value is away from the asymptotic numerical value. The GCI converts the error estimate into an error or uncertainty band, which is appropriate when one does not have a high degree of confidence in the error estimate. A small value of the GCI indicates that the computation is within the asymptotic range \cite{Roy05}.
The GCI for the fine grid solution is given as \cite{Roache94}
\begin{equation}
    \mathrm{GCI} =\frac{F_s}{g_{12}^{\hat{p}} - 1} \bigg| \frac{ V_1 - V_2 }{V_1} \bigg|,
    \label{26}
\end{equation}
where 
\begin{itemize}
    \item $F_s = 3$ if the observed order of accuracy is calculated using two different solutions.
    \item $F_s = 1.25$ if the observed order of accuracy is calculated using three different solutions and when the formal order agrees with the observed order within 10\%.
\end{itemize}
We cannot ensure that  the approximate solution via Richardson extrapolation approximates the exact solution of mathematical model well. Thus, it is beneficial to use a safety factor to avoid the failure of the designed error estimator. 
For the first case, one cannot guaranty that the solution obtained with only two grid sizes is within the asymptotic range. Thus, it is recommended to use the GCI with caution. 
In \cite{Roy10} a safety  factor of $3$ is suggested. In the second case, if the formal order of accuracy matches the observed order of accuracy, then it is recommended to use a safety factor of  $1.5$.
In this paper, when $\hat{p}$ agrees $p_f$ within 10\%, then the factor of safety of $1.25$ is used.

Algorithm \ref{Algo3} shows the steps to obtain the suitable full model dimension. The algorithm is initiated by reasonably moderate grid size $h_1^{init}$. We construct two new grid sizes $h_2,h_3$ based on the grid refinement ratios $g_{12},g_{23}$. Furthermore, a full model is solved for $h_1,h_2,h_3$ to determine the observed order of accuracy $\hat{p}$. Once we calculate $\hat{p}$, the algorithm computes the discretization error estimator $\Bar{\epsilon}_h$ and checks it against the user-defined tolerance $e_{tol}^h$. The algorithm is terminated  if $\Bar{\epsilon}_h < e_{tol}^h$ else, it refines $h_1$, and the process repeats itself until convergence or after $K$ iterations, and the last grid size is used to construct the full model.
\begin{algorithm}[\textbf{Algorithm to select the full model dimension}]
~\\
\textbf{Input:} Initial step size $h_1^{init}$, $e_{tol}^h$, $g_{12}, g_{23} > 1.3$. \\
\textbf{Output:} $h_1$.
\begin{itemize}
    \item [ ] Set $h_1 = h_1^{init}$.
    \item [ ] FOR $i=1,...,K$
        \item [ ] \quad Solve the full model with a grid size $h_1$ and store the result in $V_1$.
        \item [ ] \quad Compute $h_2 = g_{12}h_1$ and $h_3 = g_{23}g_{12}h_1$
        \item [ ] \quad Solve the full model with grid sizes $h_2$ and $h_3$ and store the results in $V_2$ and $V_3$ respectively.
        \item [ ] \quad Calculate $\hat{p}$ using Algorithm \ref{Algo2}.
        \item [ ] \quad Compute the discretization error estimator given by (\ref{25}) $\Bar{\epsilon}_h = \frac{|V_1 - V_2|}{V_1(g_{12}^{\hat{p}} - 1)}$.
        \item [ ] \quad IF $\Bar{\epsilon}_h < e_{tol}^h$
        \item [ ] \quad \quad break.
        \item [ ] \quad END.
        \item [ ] \quad Refine the mesh and choose a finer grid size $h_1$.
    \item [ ] END.
\end{itemize}
\label{Algo3}
\end{algorithm}
\subsection{Model Order Reduction Error}
\label{MORError}
For the model redeuction procedure we use a snapshot-based POD method, which represents the PDE solution in the space spanned by a reduced basis. The reduced basis is obtained by a linear combination of the snapshots generated by solving the full model for some training parameters.

Let $\mathcal{X}$ be the Hibert space given with an inner product $(\cdot, \cdot)_\mathcal{X}$ and a norm $\| \cdot\|_\mathcal{X}$. Let $V_1,V_2,\cdots,V_k \in \mathcal{X}$ be the snapshots obtained by solving the full model. We consider
\begin{equation*}
    \mathcal{V} = \mathrm{span}\{V_1, \cdots, V_k\},
\end{equation*}
as a span of the snapshots $\{V_j \}_{j=1}^{k}$, provided at least one of them is non-zero.
Let $\{\phi_i \}_{i=1}^{d}$ denote an orthonormal basis of $\mathcal{V}$ with $d= \mathrm{dim}(\mathcal{V})$. The goal of POD is to identify the most important characteristics of an ensemble of snapshots.
We can represent each (block) column of the snapshot matrix as
\begin{equation}
    V_j = \sum_{i=1}^{d} (V_j \phi_i) \phi_i, \hspace{2cm} \mathrm{for} \, j = 1,\cdots, k.
    \label{27}
\end{equation}
The POD approach obtains a reduced basis such that for every $\ell \in \{1,\cdots,d \}$, the mean square error between the snapshot $V_j$ and the corresponding $l$th partial sum of (\ref{27}) is minimized on average via
\begin{equation*}
    \min_{\{\phi_i\}_{i=1}^{\ell}} \frac{1}{k} \sum_{j=1}^{k} \| V_j - \sum_{k=1}^{\ell}(V_j \phi_i) \phi_i\|^2,
\end{equation*}
where
\begin{equation*}
    (\phi_i, \phi_j) = \delta_{ij} =
    \begin{cases}
            1, &         \text{if } i=j,\\
            0, &         \text{if } i\neq j.
    \end{cases}
    \hspace{2cm} \mathrm{for} \; 1 \leq i \leq \ell, \, 1 \leq j \leq i.
\end{equation*}
The reduced basis $\{\phi_i \}_{i=1}^d$ is obtained by computing a truncated singular value decomposition of the snapshot matrix $\hat{V}$ \cite{Golub70}
\begin{equation*}
\begin{aligned}
        \hat{V} &= \sum_{i=1}^k \Sigma_i \phi_i \psi_i,\\
        \hat{V} &= \Phi \Sigma \Psi^T.
\end{aligned}
\end{equation*}
Let $\Sigma_1 \ge \Sigma_2 \ge \cdots \Sigma_k > 0$ denote the singular values, and $\phi_1,\phi_2,\cdots,\phi_k \in \mathbb{R}^n$ the left singular vectors  of $\hat{V}$. The reduced basis $Q$ then consists of $d$ left singular vectors where the dimension $d$ is chosen such that it minimizes the projection error $\epsilon_{\mathrm{POD}}$  \cite{Kunisch01}
\begin{equation}
    \epsilon_{\mathrm{POD}} = \frac{1}{k} \sum_{j=1}^{k} \| V_j - \sum_{k=1}^{\ell}(V_j \phi_k) \phi_k\|^2 = \sum_{\ell = d+1}^{k} \Sigma_\ell^2.
    \label{28}
\end{equation}
This error characterizes the ability of the snapshot data to be represented in a low dimensional space. The goal of reduced modeling is to select the best suitable snapshots and to reduce the POD error. We use the error estimator $\epsilon_{\mathrm{POD}}$ to choose the dimension $d$ of the reduced basis in Algorithm \ref{Algo5}. Furthermore, in order to test how well the reduced model approximates the full model, we define an error between the full model and the reduced model
\begin{equation}
    \epsilon_{\mathrm{RM}} = \frac{\sqrt{\sum_{n=1}^N \|V_n - \Bar{V}_n \|_2^2}}{\sqrt{\sum_{n=1}^N \|V_n \|_2^2}},
    \label{29}
\end{equation}
where $V$ is the solution obtained using the full model and $\Bar{V}$ is the result obtained using the reduced model.

The quality of the reduced basis $Q$ depends on  $\hat{V} = [V_1(\rho_1),v_2(\rho_2),\cdots,V_k(\rho_k)]$, which depend on the training parameter set $\{\rho_1,\rho_2,\cdots, \rho_k \}$. In \cite{binder2020}, a classical and an adaptive greedy sampling approach is presented, to locate the training parameter set efficiently. We use an  a posteriori residual error estimator (\ref{5}) $\varepsilon = \|p(\cdot,\rho)\|$ to monitor the convergence of the classical greedy sampling algorithm and to terminate the training procedure. In the adaptive greedy sampling approach, the algorithm constructs an error model $\Bar{\epsilon}_{\mathrm{RM}}$ based on the available data $E_p$ from the residual errors $\{\varepsilon_i^{\mathrm{bef},\mathrm{aft}}\}_{i=1}^{I_{max}}$ and the reduced model errors $\{{\epsilon^{\mathrm{bef},\mathrm{aft}}_{\mathrm{RM},i}}\}_{i=1}^{I_{max}}$ and uses this error model to monitor the convergence of the algorithm. Here the superscripts bef and aft show the error values obtained before and after updating the reduced basis.

For this approach, one full model is solved for an optimal parameter group $\rho_I$ at each iteration $i$ of the classical and the adaptive greedy procedures. The solution obtained for $\rho_I$ is then stored as $V_i(\rho_I)$ and updates the snapshot matrix $\hat{V} = [V_\ell(\rho_I)]_{\ell=1}^{i}$
\begin{equation*}
    \hat{V}_\ell = [V_1(\rho_I^1), V_2(\rho_I^2), \cdots, V_\ell (\rho_I^\ell)] \in \mathbb{R}^{M\times (T \times i)} \hspace{2cm} \ell = 1,...,i,
\end{equation*}
where $T$ is the maturity. Now we consider
\begin{equation*}
    \mathcal{V}_\ell = \mathrm{span}\{V_1(\rho_I^1), V_2(\rho_I^2), \cdots, V_\ell(\rho_I^\ell)\},
\end{equation*}
as a span of these snapshots and obtain the reduced basis by computing its truncated singular value decomposition. At each greedy iteration, a new reduced basis $Q_i$ is obtained of dimension $d$ and the projection error based on (\ref{28}) is
\begin{equation}
    \epsilon_{\mathrm{POD}}^{\mathrm{AG,CG}} = \frac{1}{i\times T} \sum_{i=1}^{i \times T} \| V_i(\rho_I) - \sum_{i\times T=1}^{\ell}(V_i(\rho_I) \phi_k) \phi_k\|_2^2 = \sum_{\ell = d+1}^{i\times T} \Sigma_\ell^2,
    \label{30}
\end{equation}
which is the sum of squares of neglected singular values.

The computation of the SVD of an $m\times n$ matrix 
has an overall computational complexity of  $\mathcal{O}(mn^2)$ \cite{GolV13}, so for 
the matrix $\hat{V} \in \mathbb{R}^{M \times (T \times i)}$ it is
\begin{equation*}
    \mathcal{O}(M(T\times i)^2) = \mathcal{O}(MT^2i^2).
\end{equation*}
After $i$ iterations, the computational complexity due to the series of SVDs for the snapshot matrix will be
\begin{equation*}
    \mathcal{O}(\sum_{\ell=1}^i M \times (T \times i)^2 ) = \mathcal{O}(MT^2 i^3).
\end{equation*}
We see that the complexity increases substantially with each proceeding iteration, which motivates us to implement a truncated SVD approach. Recent studies have shown that the randomized algorithms in some applications  \cite{Feng18,Halko11} yield an incredible speed-up. It is proposed to compute a truncated singular value decomposition using randomized algorithms in the following two steps. The approach starts by computing an approximate basis $G \in \mathbb{R}^{M\times k}$ with $k$ orthonormal columns for the range of $\hat{V}$ such that
\begin{equation*}
\begin{aligned}
        \hat{V} \approx G G^* \hat{V},\\
        \|\hat{V} - GG^* \hat{V} \| < \epsilon_{tol},
\end{aligned}
\end{equation*}
where $k < T\times i$ and $\epsilon_{tol}$ is some specified tolerance. In the second step  the matrix $G$ is used to compute the truncated SVD of the matrix $\hat{V}$. Let $\Tilde{V} = G^* \hat{V} \in \mathbb{R}^{k \times (T\times i)}$ be the reduced matrix determined  from$\hat{V}$. Now we compute the SVD of the reduced matrix
\begin{equation*}
    \Tilde{V} = G^* \hat{V} = \Tilde{\Phi} \Tilde{\Sigma} \Tilde{\Psi}^T
\end{equation*}
which leads to a substantial speedup.
Multiplying both sides by $G$, we get
\begin{equation}
    \begin{aligned}
        G G^* \hat{V} = G \Tilde{\Phi} \Tilde{\Sigma} \Tilde{\Psi}^T,\\
        \hat{V} \approx G \Tilde{\Phi} \Tilde{\Sigma} \Tilde{\Psi}^T.
    \end{aligned}
    \label{31}
\end{equation}
Considering $\Phi = G \Tilde{\Phi}$, we get our left singular vector matrix, which we can use to obtain a reduced basis. 

We can achieve even more speedup by incorporating the economy-sized (compact) SVD of the matrix $\Tilde{V}$, which removes extra rows or columns of zeros from the diagonal matrix of singular values and the columns in either $\Tilde{\Phi}$ or $\Tilde{\Psi}$. Removing these zeros and columns can improve execution time and reduces storage requirements without compromising the accuracy of the decomposition.
The SVD determines the rank $k$ of a matrix, as the number of nonzero singular values of a matrix $\Tilde{V} \in \mathbb{R}^{k\times {T\times i}}$. 
Then its compact SVD reads
\begin{equation*}
    \Tilde{V} = \Tilde{\Phi} \Tilde{\Sigma} \Tilde{\Psi}^T,
\end{equation*}
The left and right singular vector matrices are $\Tilde{\Phi} \in \mathbb{R}^{k\times k}$, and $\Tilde{\Psi} \in \mathbb{R}^{(T\times i) \times k}$. Thus, the desired left singular matrix is $\Phi = G\Tilde{\Phi} \in \mathbb{R}^{M \times k}$, which we use to obtain the reduced basis $Q$.
We obtain the orthonormal matrix $G$ with as few columns as possible such that
\begin{equation*}
    \|\hat{V} - GG^* \hat{V} \| < \varepsilon_{tol}
\end{equation*}
where $\varepsilon_{tol}$ is some tolerance and $\|\cdot\|$ is the $L_2$-norm. The range of $G$ is a $k$ dimensional subspace that represents most of the action of $\hat{V}$, and we select $k$ as small as possible. We can use the SVD to decide the dimension $k$. Let $\Sigma_j$ be the $j$th singular value of the matrix $\hat{V}$. For each $j \geq 0$
\begin{equation*}
    \underset{\mathrm{rank}(GG^* \hat{V}) \leq j}{\mathrm{min}} \| \hat{V} - G G^* \hat{V} \| = \Sigma_{j+1},
\end{equation*}
where the columns of $G$ are $k$ dominant singular vectors of $\hat{V}$. The number of singular values of $\hat{V}$ that exceeds the tolerance $\varepsilon_{tol}$ makes the minimal rank $k$. Draw $k$ Gaussian random vectors $\hat{G} = [\hat{g}_1, \cdots, \hat{g}_k] \in \mathbb{R}^{(T\times i)\times k}$ and project them on to the linear map of a matrix $\hat{V}$ such that
\begin{equation*}
    \begin{aligned}
        y_\ell &= \hat{V} \hat{g}_\ell \hspace{2cm} \ell = 1,...,k,\\
        Y &= \hat{V}\hat{G}.
    \end{aligned}
\end{equation*}
We find the orthonormal basis of the matrix $Y \in \mathbb{R}^{M \times k}$, which is nothing but our desired matrix $G$,  using a Gram–Schmidt like algorithm or QR factorization.
We can determine $G$ adaptively using the error $\|(I - GG^*) \hat{V}\|$. We can obtain some information about this error by calculating $\|(I - GG^*) \hat{V}\hat{g}\|$. For $k$ Gaussian random vectors, \cite{Halko11} gives an posterior error estimator as
\begin{equation}
    \|\hat{V} - GG^* \hat{V}\| \leq 10\sqrt{\frac{2}{\pi}} \underset{i=1,\dots,k}{\mathrm{max}}\|(I - GG^*) \hat{V}\hat{g}_i\|.
    \label{32}
\end{equation}
This error estimate can be combined with any method for constructing an approximate basis for the range of $Y$, see \cite{Halko11} for more details. The adaptive algorithm \ref{Algo4} combined with Gram-Schmidt orthonormalization is given below. The algorithm is initiated by drawing $k$ Gaussian random vectors of length $T\times i$ and computes the orthonormal set of vectors $\{g_1,...,g_k\}$ that form the matrix $G$. We stop the procedure once the largest value in a set $\{\|y_i\|\}_{i}^{i+k}$ is less than $\epsilon/(10\sqrt{2/\pi})$.
\begin{algorithm}[\textbf{Algorithm to obtain the orthonormal matrix $G$}]
~\\
\textbf{Input:} Snapshot matrix $\hat{V}$. \\
\textbf{Output:} $G$.\\
\textbf{Function:} $GSorth()$
\begin{itemize}
    \item [ ] Draw Gaussian random vectors $\{\hat{g}_i\}_{i=1}^{k}$.
    \item [ ] FOR $i=1,\dots,k$
        \item [ ] \quad $y_i = \hat{V} \hat{g}_i$
        \item [ ] \quad FOR $j = 1,\dots,i-1$
        \item [ ] \quad \quad Compute $y_i = y_i - (g_j^T\hat{g}_i) g_j$.
        \item [ ] \quad END
        \item [ ] \quad $g_i = \frac{y_i}{\|y_i\|}$.
        \item [ ] \quad Set $G_i = [g_i]$.
        \item [ ] \quad Draw a new Gaussian random vector $\hat{g}_{i + k}$.
        \item [ ] \quad Compute $y_{i+k} = (I - G_iG_i^T)\hat{V}\hat{g}_{i+k}$.
        \item [ ] \quad IF $\mathrm{max}\{\|y_i\|,\|y_{i+1}\| \dots, \|y_{i+k}\|\} < \epsilon/(10\sqrt{2/\pi})$
        \item [ ] \quad \quad STOP.
        \item [ ] \quad END IF
    \item [ ] END
    \item [ ] $G = G_i$.
\end{itemize}
\label{Algo4}
\end{algorithm}
We use the randomized singular value decomposition in the greedy sampling approach to determine the reduced basis, as shown in Algorithm \ref{Algo5}. The algorithm constructs a snapshot matrix $\hat{V}$ composed of full model solutions for each optimal parameter group at the $i$th iteration. It then selects $k$ Gaussian random vectors and obtains the matrix $G$, using the function $GSorth()$.
\begin{algorithm}[\textbf{Algorithm to obtain a reduced basis using a randomized singular value decomposition}]
~\\
\textbf{Input:} Snapshot matrix $\hat{V}$, $e_{tol}^d$. \\
\textbf{Output:} $Q$.
\begin{itemize}
    \item [ ] Construct a snapshot matrix $\hat{V} = \{V(\rho_\ell) \}_{\ell}^i$ at $i$th iteration.
    \item [ ] Obtain $G$ using 'GSorth' function defined in Algorithm \ref{Algo4}.
    \item [ ] Compute $\Tilde{V} = G^T \hat{V}$.
    \item [ ] Compute a compact singular value decomposition of $\Tilde{V} = \Tilde{\Phi} \Tilde{\Sigma} \Tilde{\Psi}$.
    \item [ ] Set $\Phi = G \Tilde{\Phi}$.
    \item [ ] Choose a parameter group $\rho_{\mathrm{test}}$ for testing.
    \item [ ] Solve the full model for $\rho_{\mathrm{test}}$ and store result in $V$.
    \item [ ] $\Sigma = [\Sigma_1,...,\Sigma_k]$.
    \item [ ] FOR $j=1,...,\mathrm{length}(\mathrm{diag}(\Sigma))$
        \item [ ] \quad Compute $\epsilon_{\mathrm{POD}} = \sum_{\ell = j+1}^{k} \Sigma_\ell^2$.
        \item [ ] \quad Solve the reduced model for $\rho_{\mathrm{test}}$ for store result in $\Bar{V}$.
        \item [ ] \quad Compute
    $\epsilon_{\mathrm{RM}} = \frac{\sqrt{\sum_{n=1}^N \|V_n - \Bar{V}_n \|_2^2}}{\sqrt{\sum_{n=1}^N \|V_n \|_2^2}}$
        \item [ ] \quad Compute $\epsilon_d = \epsilon_{\mathrm{RM}} + \epsilon_{\mathrm{POD}}$.
        \item [ ] \quad IF $\epsilon_{d} < e_{tol}^d$
        \item [ ] \quad \quad d = j
        \item [ ] \quad \quad STOP
        \item [ ] \quad END IF
    \item [ ] END.
    \item [ ] $Q = [\phi_1 ,..., \phi_d]$.
\end{itemize}
\label{Algo5}
\end{algorithm}
A low dimensional matrix $\Tilde{V}$ is constructed using the newly obtained $G$ with only $k$ columns. The algorithm computes a compact SVD of $\Tilde{V}$ and sets $\Phi = G\Tilde{\Phi}$ to obtain a reduced basis $Q$. Finally, the dimension of a reduced basis is selected that minimizes the projection error and the reduced model error.
\subsection{Sampling Error}
\label{SampError}
To design a surrogate model to locate the training parameters we will use principal component regression. There are two errors associated with the principal component regression model, i.e., a projection error due to the principal component analysis (PCA) and an error of prediction for the surrogate model. 

Let us first discuss the projection error. PCA reduces the problem of estimating $m$ coefficients to the more straightforward problem of determining $p$ coefficients by selecting a few relevant principal components. It is a technique for projecting multidimensional data onto a low dimensional subspace with minimal loss of variance. We determine the projection error associated with the PCA in a similar sense to the POD.
At the $k$th iteration we have a data matrix $\hat{\mathcal{P}}_k = [\rho_1,\dots,\rho_{c_k}]$ composed of $c_k$ parameter groups. To build the surrogate model, we compute an SVD of the matrix $\hat{\mathcal{P}}_k$
\begin{equation*}
    \hat{\mathcal{P}}_k = \hat{\Phi} \hat{\Sigma} \hat{\Psi},
\end{equation*}
where the columns of $\hat{\mathcal{P}}_k\hat{\Psi}$ are principal components. We use only $p$ principal components to construct a fairly accurate surrogate model.
Similar to (\ref{28}), the projection error associated to PCA, if only $p$ principal components are used, is
\begin{equation}
    \epsilon_{\mathrm{PCA}} = \frac{1}{m} \sum_{j=1}^{m} \| \rho_j - \sum_{m=1}^{\ell}(\rho_j \phi_k) \phi_k\|^2 = \sum_{\ell = p+1}^{m} \hat{\Sigma}^2_\ell.
    \label{33}
\end{equation}
In order to obtain a satisfactory regression model, we aim to reduce the projection error and select the principal components accordingly.

In the following, we discuss the mean squared error of prediction, which estimates the quality of the regression model, i.e., the surrogate model. 
The error estimators for surrogate models are analyzed onindependent test data. It may happen that the test set is not large enough, as in our case, during the start of the greedy iteration.
If $c_k$ is not much larger than $m$, then the model may give weak predictions due to the risk of over-fitting for the parameter groups which are not used in model training.
In such situations, the error estimators are designed on the learning data where leave-one-out cross-validation is an elegant choice.
As we have a data set with $c_k$ observations, the leave-one-out cross-validation approach trains the model on $c_k - 1$ data points, and test data contains one observation. We repeat this process $c_k$ times and compute the mean square error. The major advantage of this approach is that it has far less bias as if we use the entire data for training.
However, the computational cost of one full leave-one-out cross-validation is high, and it may produce inconsistent and varying results \cite{Wong15}. Additionally, it tends to have high variance and generates very different estimates if repeated estimates with different initial samples of data from the same distribution.

K-fold cross-validation is an alternative advocated to circumvent these drawbacks \cite{Fushiki11}. This technique involves randomly dividing the data set into $K$ subsets known as \emph{folds} of approximately equal size. One fold is kept for testing, and the model is trained on the remaining $K-1$ folds. The process repeats for $K$ iterations each time with $K-1$ folds as the training set and the $K$th fold as a test set. As we repeat the process $K$ times, we get $K$ errors per each fold, which are then used to compute the mean square error of prediction. The K-fold cross-validation is computationally more efficient than the leave-one-out cross-validation as we repeat the process only for $K$ iterations instead of $c_k$ iterations.

Let $L = [\hat{\mathcal{P}}_k, \hat{\varepsilon}] \in \mathbb{R}^{c_k \times (m+1)}$ be a learning data matrix consisting of $c_k$ parameter groups and error estimator values at $k$th iteration of the adaptive greedy algorithm. Here $\hat{\mathcal{P}}_k = [\rho_1, \cdots, \rho_{c_k}] \in \mathbb{R}^{c_k \times m}$ and $\hat{\varepsilon} = [\varepsilon_1, \cdots, \varepsilon_{c_k}] \in \mathbb{R}^{c_k \times 1}$. Let the dataset $L$ be partitioned into $K$ groups or folds $L_1, \cdots, L_K$ such that $L_i \cap L_j =  \emptyset$ for any $i \neq j$. Without loss of generality, we suppose that $c_k$ is a multiple of $K$ where $c_k$ denotes the size of dataset $L$. The size of each fold is then $c_k/K$. We train our surrogate model $\Bar{\varepsilon}$ on $L\setminus L_K$, i.e., all data available in the matrix $L$ except in $L_K$, while the data set $L_K$ is used as a test set. After $K$ iterations, we take the sum of errors and calculate its mean to get the $\mathrm{MSEP}_{K-\mathrm{CV}}$. The K-fold cross-validation error is then given as \cite{mevik04}
\begin{equation}
    \mathrm{MSEP}_{K-\mathrm{CV}} = \frac{1}{c_k} \sum_{k =1}^{K} \sum_{\substack{i \in L_K\\ i = 1}}^{c_k/K} (\Bar{\varepsilon}(\rho_i) - \varepsilon_i)^2,
    \label{34}
\end{equation}
where the inner sum is calculated over the $K$th fold. The bias of $\mathrm{MSEP}_{K-\mathrm{CV}}$ is $\mathcal{O}(K-1)^{-1} c_k^{-1}$ \cite{Fushiki11}.
In K-fold cross-validation, the surrogate model is trained on $L\setminus L_K$ data instead of all $L$. It can be expected that this surrogate model may be inferior as compared to the model trained on all data $L$. Thus, it is suggested to use the adjusted K-fold cross-validation to avoid this problem. The adjustment is \cite{mevik04}
\begin{equation}
    \mathrm{MSEP}_{adj.} = \mathrm{MSEP}_{app.} - \frac{1}{c_k} \sum_{k=1}^K \frac{1}{K} \sum_{\substack{i \not\in L_K \\ i=1}}^{c_k - c_k/K} (\Bar{\varepsilon}(\rho_i) - \varepsilon_i)^2,
    \label{35}
\end{equation}
where $\mathrm{MSEP}_{app.}$ is called the apparent MSEP or mean squared error of calibration (MSEC), which uses the whole learning data $L$ as a test set. The second term of (\ref{35}) is calculated as follows. In contrast to (\ref{34}), instead of testing the model trained on the excluded data set $L_k$, we test it on the training data set $L\setminus L_k$ itself. After $K$ iterations, we compute the mean by dividing the sum of $K$ errors by $c_k$, which is nothing but the second term of (\ref{35}).
The adjustment is then the difference between the apparent error and the error calculated by training and testing the model on the $L\setminus L_K$ data. In short, this adjustment considers the errors associated with the K folds, which are not used for the training purpose during the K-fold cross-validation. The apparent MSEP is given as
\begin{equation}
    \mathrm{MSEP}_{app.} = \frac{1}{c_k} \sum_{i=1}^{c_k} (\Bar{\varepsilon}(\rho_i) - \varepsilon_i)^2.
    \label{36}
\end{equation}
Thus, the adjusted K-fold cross validation error estimate is
\begin{equation}
\begin{aligned}
         &\mathrm{MSEP}_{K-\mathrm{CV}adj.} = \mathrm{MSEP}_{K-\mathrm{CV}} + \mathrm{MSEP}_{adj.},\\
         &= \frac{1}{c_k} \sum_{k =1}^{K} \sum_{\substack{i \in L_K\\ i = 1}}^{c_k/K} (\Bar{\varepsilon}(\rho_i) - \varepsilon_i)^2 + \frac{1}{c_k} \sum_{i=1}^{c_k} (\Bar{\varepsilon}(\rho_i) - \varepsilon_i)^2 - \frac{1}{c_k} \sum_{k=1}^K \frac{1}{K} \sum_{\substack{i \not\in L_K \\ i=1}}^{c_k - c_k/K} (\Bar{\varepsilon}(\rho_i) - \varepsilon_i)^2.
\end{aligned}
\label{37}
\end{equation}
We can now define the total sampling error associated with the PCR technique as a summation of the projection error and the adjusted K-fold cross-validation error
\begin{equation}
    \epsilon_{samp} = \epsilon_{\mathrm{PCA}} + \mathrm{MSEP}_{K-\mathrm{CV}adj.}.
    \label{38}
\end{equation}
Algorithm \ref{Algo6} presents a methodology to obtain a suitable surrogate model that minimizes the sampling error. It starts by computing a projection error $\epsilon_{\mathrm{PCA}}$ associated with the PCA. The surrogate model $\Bar{\varepsilon}$ is then constructed using $j=p$ principal components. To assess the quality of the surrogate model, we determine the cross-validation error. The data set $L$ is divided into $K$ folds. We train the designed surrogate model on all data available in the matrix $L$ except in $L_k$. Furthermore, the algorithm evaluates the surrogate model on the omitted data $L_k$, then computes errors, and subsequently determines the K-fold cross-validation error. $\mathrm{MSEP}_{K-\mathrm{CV}}$ is then adjusted and added to $\epsilon_{\mathrm{PCA}}$ to get the sampling error $\epsilon_{samp}$. If the sampling error is less than the user-defined tolerance $e_{tol}^{samp}$, then the algorithm terminates. Finally, we use the surrogate model built using $j=p$ principal components as explained in \cite{binder2020}.
\begin{algorithm}[\textbf{Algorithm to estimate the sampling error}]
~\\
\textbf{Input:} $L = [\hat{P}_k, \hat{\varepsilon}]$, $e_{tol}^{samp}$, K. \\
\textbf{Output:} $\epsilon_{samp}$ (\ref{38}).
\begin{itemize}
    \item [ ] Compute an SVD of the matrix $\hat{\mathcal{P}}_k = \hat{\Phi} \hat{\Sigma} \hat{\Psi}$
    \item [ ] Let $\hat{\Sigma} = [\hat{\Sigma}_1,\dots,\hat{\Sigma}_k]$.
    \item [ ] FOR $j=1,...,\mathrm{length}(\mathrm{diag}(\hat{\Sigma}))$
        \item [ ] \quad Compute $\epsilon_{\mathrm{PCA}} = \sum_{\ell = j+1}^{k} \hat{\Sigma}_\ell^2$.
        \item [ ] \quad Build a surrogate model $\Bar{\varepsilon}$ using $p=j$ principal components.
        \item [ ] \quad Divide $L$ into $K$ random folds.
        \item [ ] \quad FOR $k = 1,...,K$
        \item [ ] \quad \quad Train the surrogate model on $L\setminus L_k$ data.
        \item [ ] \quad \quad Solve the trained model on $L_k$ data set.
        \item [ ] \quad \quad Compute squared errors $\sum_{\substack{i \in L_K\\ i = 1}}^{c_k/K} (\Bar{\varepsilon}(\rho_i) - \varepsilon_i)^2$.
        \item [ ] \quad END
        \item [ ] \quad Compute (\ref{34}) $\mathrm{MSEP}_{K-\mathrm{CV}}$.
        \item [ ] \quad Compute (\ref{36}) $\mathrm{MSEP}_{app}$.
        \item [ ] \quad Compute (\ref{35}) $\mathrm{MSEP}_{adj}$.
        \item [ ] \quad Compute (\ref{37}) $\mathrm{MSEP}_{K-\mathrm{CV}adj.}$
        \item [ ] \quad Compute (\ref{38}) $\epsilon_{samp} = \epsilon_{\mathrm{PCA}} + \mathrm{MSEP}_{K-\mathrm{CV}adj.}$
        \item [ ] \quad IF $\epsilon_{samp} < e_{tol}^{samp}$
        \item [ ] \quad \quad p = j
        \item [ ] \quad \quad STOP
        \item [ ] \quad END IF
    \item [ ] END.
\end{itemize}
\label{Algo6}
\end{algorithm}
\subsection{Sensitivity Analysis}
\label{Sensitivity}
In the previous subsections, we addressed the errors $\epsilon_h, \epsilon_{\mathrm{POD}}, \epsilon_{\mathrm{RM}},\epsilon_{samp}$ associated with each stage of the model order reduction approach. This section presents a sensitivity analysis for these errors that provides quantitative information regarding the relative contribution of each error to the model output. Sensitivity analysis can be helpful in various situations, including forecasting as well as identifying where improvements or adjustments are needed in a process \cite{Frey02}. Such information helps to improve the model output by selecting the appropriate grid size $h$, the reduced dimension $d$, and the optimal sampling parameters $c$.
The sensitivity indices rank the numerical errors according to their contribution to the total error. We can use this ranking to allocate the computational resource efficiently. More computational resources could be provided to the most sensitive parameters \cite{Liang11}.

In general, there are two types of sensitivity analysis, i.e., local and global \cite{Borgonovo16,Frey02}. The local sensitivity analysis evaluates changes in the output with respect to the variations in a single factor by keeping other factors constant \cite{Iooss15}. Such sensitivity is often evaluated through gradients or partial derivatives of the output functions at these factors. The main limitation of the local sensitivity analysis is that it evaluates the factors one at a time and does not consider the interactions and simultaneous changes in all factors \cite{Saltelli10}.
Global sensitivity analysis overcomes this drawback by quantifying the importance of all factors and their interactions to the model output. It provides an overall influence of factors on the output in contrast to a local view of partial derivatives as in local sensitivity analysis.

It is easy and computationally feasible to compute the local sensitivity of deterministic errors such as the discretization error by computing the norm of the first-order factor in a Taylor series expansion or the partial derivative using finite differences \cite{Frey02}. For stochastic errors like the sampling error, it is more appropriate to compute the change in the output with respect to the variations in inputs \cite{Liang11}. However, this creates a difficulty in comparing the relative contributions of the deterministic errors versus the stochastic errors. Thus, it is necessary to propose an approach that can be applied to both types of errors.

There are different local and global sensitivity analysis methods, such as the one-at-a-time (OAT) approach \cite{Borgonovo16}, the automatic differentiation technique (AD) \cite{Corliss02}, variance-based sensitivity analysis (the Sobol method) \cite{Mara82,Saltelli10}, and Fourier amplitude sensitivity analysis (FAST) \cite{McRae82}.
The local sensitivity analysis approaches like OAT and AD calculate the first-order partial derivatives of the outputs with respect to small changes in the inputs. In AD, a computationally efficient computer code automatically evaluates the partial derivatives and can be applied without having detailed knowledge of the algorithm implemented in the model. However, the technique may be limited to specific computer languages and requires specific libraries, and the accuracy for sensitivity results depends on the numerical method used in the AD software. Both the Sobol method and the FAST method  are variance-based sensitivity analysis techniques capable of computing first-order (local) and global sensitivity indices. The main difference between these two approaches is the underlying algorithm in multidimensional integration of sensitivity indices. The FAST approach uses a pattern search technique based on a sinusoidal function, while the Sobol approach uses a Monte Carlo technique to determine the sensitivities. The disadvantage of the FAST approach is that it demands specific sampling and may perform poorly for discrete inputs and models with discontinuities.
Among the discussed techniques, the Sobol method is one of the most powerful techniques being model-independent and robust, and it considers interactions between different input factors. However, it may become computationally intensive if there are many input factors as it uses Monte Carlo techniques. In our case, the number of input factors (numerical errors) is small, which eases the computational burden.
Moreover, \cite{Liang11} showed that the variance-based sensitivity analysis approach could be applied considering both deterministic and stochastic numerical errors.
Thus, in this work, we use the Sobol method to perform the sensitivity analysis.

Consider a model output $Y$ depending on some input factors $(X_1, \dots ,X_p)$
\begin{equation*}
    Y = f(X_1,\dots,X_p).
\end{equation*}
Now, the variance based first order effect of $X_i$ is \cite{Mara82}
\begin{equation}
    \mathrm{Var}_{X_i} \Big( \mathbb{E}_{X_{\sim i}} (Y | X_i) \Big), \; \; \; \; i = 1,\dots,p,
    \label{39}
\end{equation}
where $X_i$ is the $i$th factor and $X_{\sim i}$ is the $p-1$ dimensional space composed of all factors except $X_i$. The inner expectation operator means that the mean of $Y$ is taken over all possible values of $X_{\sim i}$ while the factor $X_i$ is kept fixed. $\mathrm{Var}_{X_i}(\cdot)$ shows that the variance is taken over all possible values of $X_i$. Using the law of total variance, we can write \cite{Mood74}
\begin{equation}
    \mathrm{Var}(Y) = \mathrm{Var}_{X_i} \Big( \mathbb{E}_{X_{\sim i}} (Y | X_i) \Big) + \mathbb{E}_{X_i} \Big( \mathrm{Var}_{X_{\sim i}} (Y | X_i) \Big),
    \label{40}
\end{equation}
by normalizing
\begin{equation}
    1 = \frac{\mathrm{Var}_{X_i} \Big( \mathbb{E}_{X_{\sim i}} (Y | X_i) \Big)}{\mathrm{Var}(Y)} + \frac{\mathbb{E}_{X_i} \Big( \mathrm{Var}_{X_{\sim i}} (Y | X_i) \Big)}{\mathrm{Var}(Y)},
    \label{41}
\end{equation}
where the first order (local) sensitivity index for the factor $X_i$ is given by the first term of (\ref{41})
\begin{equation}
    S_i =  \frac{\mathrm{Var}_{X_i} \Big( \mathbb{E}_{X_{\sim i}} (Y | X_i) \Big)}{\mathrm{Var}(Y)}.
    \label{42}
\end{equation}
Here $S_i$ is a normalized index as $\mathrm{Var}_{X_i} \Big( \mathbb{E}_{X_{\sim i}} (Y | X_i) \Big)$ varies from zero to $\mathrm{Var}(Y)$ and measures the first order effect of $X_i$ on the output $Y$. $\mathbb{E}_{X_i} \Big( \mathrm{Var}_{X_{\sim i}} (Y | X_i) \Big)$ is the residual \cite{Saltelli10}.
It is necessary to note that the low value of the first-order sensitivity index does not mean that the factor is not important. It might make a more significant contribution to the output by interactions with other factors. Thus, the total sensitivity index is necessary, which considers the main effect of the factor and the effect of its interactions with other factors on the output. In the sensitivity analysis framework, functional decomposition of the total variance $\mathrm{Var}(Y)$ is referred to as functional analysis of variance (ANOVA) and given as \cite{Sobol01}
\begin{equation}
    \mathrm{Var}(Y) = \sum_{i = 1}^{p} \mathrm{Var}_i(Y) + \sum_{i < j}^p \mathrm{Var}_{ij}(Y) + \cdots + \mathrm{Var}_{1...p}(Y),
    \label{43}
\end{equation}
where
\begin{equation*}
\begin{aligned}
    \mathrm{Var}_i(Y) &= \mathrm{Var}_{X_i} \Big( \mathbb{E}_{X_{\sim i}} (Y | X_i) \Big),\\
    \mathrm{Var}_{ij}(Y) &= \mathrm{Var}_{X_i,X_j} \Big( \mathbb{E}_{X_{\sim i,j}} (Y | X_i,X_j) \Big) - \mathrm{Var}_i(Y) - \mathrm{Var}_j(Y),\\
    \mathrm{Var}_{ijk}(Y) &= \mathrm{Var}_{X_i,X_j,X_k} \Big( \mathbb{E}_{X_{\sim i,j,k}} (Y | X_i,X_j,X_k) \Big) - \mathrm{Var}_i(Y) - \mathrm{Var}_j(Y) - \mathrm{Var}_k(Y)\\
    &- \mathrm{Var}_{ij}(Y) - \mathrm{Var}_{ik}(Y) - \mathrm{Var}_{jk}(Y),\\
    \vdots\\
    \mathrm{Var}_{1,...,p} &= \mathrm{Var}(Y) - \sum_{i=1}^p \mathrm{Var}_i(Y) - \sum_{1\leq i < j \leq p} \mathrm{Var}_{ij}(Y) - \cdots - \sum_{1\leq i1<...<ip-1\leq p} \mathrm{Var}_{i1,...,ip-1}.
\end{aligned}
\end{equation*}
The first-order sensitivity index (\ref{42}) can then be obtained from the first $p$ terms of the decomposition (\ref{43}) as
\begin{equation*}
    S_i = \frac{\mathrm{Var}_i(Y)}{\mathrm{Var}(Y)}.
\end{equation*}
The other terms of the decomposition can be interpreted as the higher-order sensitivity indices. The second-order index $S_{ij}$ represents the effect of interactions of factors $X_i$ and $X_j$ on the output
\begin{equation*}
    S_{ij} = \frac{\mathrm{Var}_{ij}(Y)}{\mathrm{Var}(Y)},
\end{equation*}
and so on up to  order $p$. Therefore, for $p$ factors, we have $2^p -1$ sensitivity indexes. Based on the decomposition (\ref{43}), \cite{Homma96,Saltelli02} introduced the so-called total indices or total effects as
\begin{equation*}
    S_{T_i} = S_i + \sum_{j\neq i} S_{ij} + \sum_{j\neq i, k \neq i, j<k} S_{ijk} + \cdots = \sum_{\ell \# i} S_\ell,
\end{equation*}
where $\# i$ denotes all the indices associated with the factor $X_i$. The total sensitivity is derived using the expressions in (\ref{43}) as \cite{Saltelli10}
\begin{equation*}
    S_{T_i} = 1 - \frac{\mathrm{Var}_{X_{\sim i}} (\mathbb{E}_{X_i} (Y|X_{\sim i}))}{\mathrm{Var}(Y)}.
\end{equation*}
Using the law of total variance, similar to (\ref{40}), we can write the total variance of $Y$ by exchanging $X_i$ and $X_{\sim i}$ as
\begin{equation}
    \mathrm{Var}(Y) = \mathrm{Var}_{X_{\sim i}} \Big( \mathbb{E}_{X_{i}} (Y | X_{\sim i}) \Big) + \mathbb{E}_{X_{\sim i}} \Big( \mathrm{Var}_{X_{i}} (Y | X_{\sim i}) \Big).
    \label{44}
\end{equation}
Dividing (\ref{44}) by $\mathrm{Var}(Y)$
\begin{equation*}
    1 = \frac{\mathrm{Var}_{X_{\sim i}} \Big( \mathbb{E}_{X_{i}} (Y | X_{\sim i}) \Big)}{\mathrm{Var}(Y)} + \frac{\mathbb{E}_{X_{\sim i}} \Big( \mathrm{Var}_{X_{i}} (Y | X_{\sim i}) \Big)}{\mathrm{Var}(Y)},
\end{equation*}
we obtain
\begin{equation*}
    S_{T_i} = 1 - \frac{\mathrm{Var}_{X_{\sim i}} (\mathbb{E}_{X_i} (Y|X_{\sim i}))}{\mathrm{Var}(Y)} = \frac{\mathbb{E}_{X_{\sim i}} \Big( \mathrm{Var}_{X_{i}} (Y | X_{\sim i}) \Big)}{\mathrm{Var}(Y)}.
\end{equation*}
\subsubsection{Calculation of Sobol Indices}
The methodology to determine the local and global Sobol indices with a single simulation is presented in \cite{Saltelli10}. Consider two independent sampling matrices $\Bar{X}$ and $\hat{X}$ to compute the values of the output $Y$ corresponding to different input factors $X_1,...,X_p$.
\begin{equation*}
\Bar{X} =
\begin{bmatrix}
\Bar{X}_{11} & \cdots & \Bar{X}_{1p}\\
\vdots & \vdots & \vdots\\
\Bar{X}_{n1} & \cdots & \Bar{X}_{np}\\
\end{bmatrix}_{n \times p},
\hat{X} =
\begin{bmatrix}
\hat{X}_{11} & \cdots & \hat{X}_{1p}\\
\vdots & \vdots & \vdots\\
\hat{X}_{n1} & \cdots & \hat{X}_{np}\\
\end{bmatrix}_{n \times p},
\end{equation*}
where $n$ is the sample size and the $p$ is the number of input factors. Then each row of the matrices $\Bar{X}$ and $\hat{X}$ is then a one parameter sample set. We introduce a new matrix $\Bar{X}^i_{\hat{X}}$ whose  columns are from the matrix $\Bar{X}$ expect the $i$th column, which is from the matrix $\hat{X}$.
\begin{equation*}
\Bar{X}^i_{\hat{X}} =
\begin{bmatrix}
\Bar{X}_{11} & \Bar{X}_{12} & \cdots & \hat{X}_{1i} & \cdots & \Bar{X}_{1p}\\
\vdots & \vdots & \vdots & \vdots & \vdots & \vdots\\
\vdots & \vdots & \vdots & \vdots & \vdots & \vdots\\
\Bar{X}_{n1} & \Bar{X}_{n2} & \cdots & \hat{X}_{ni} & \cdots & \Bar{X}_{np}\\
\end{bmatrix}_{n \times p},
\end{equation*}
Based on these matrices $\Bar{X},\hat{X}$, and $\Bar{X}^i_{\hat{X}}$, we determine the local sensitivity and the global sensitivity indices as follows. We obtain a variance of $Y$ as in \cite{Mood74}
\begin{equation}
    \mathrm{Var}(Y) = \mathbb{E}(Y^2) + \mathbb{E}^2(Y).
    \label{45}
\end{equation}
Applying (\ref{45}) to $\mathrm{Var}_{X_i} \Big( \mathbb{E}_{X_{\sim i}} (Y | X_i) \Big)$, we get
\begin{equation*}
    \mathrm{Var}_{X_i} \Big( \mathbb{E}_{X_{\sim i}} (Y | X_i) \Big) = \int \mathbb{E}^2_{X_{\sim i}} (Y|X_i) dX_i - \bigg( \int \mathbb{E}_{X_{\sim i}} (Y|X_i) dX_i\bigg)^2 .
\end{equation*}
From \cite{Saltelli2002} we then have
\begin{equation}
\begin{aligned}
        \int \mathbb{E}^2_{X_{\sim i}} (Y|X_i) dX_i &= \int \int f(\hat{X}_1,\dots,\hat{X}_p) \times f(\Bar{X}_1,\dots,\hat{X}_i,\dots,\Bar{X}_p) d\hat{X} d\Bar{X}_{\sim i},\\
        &= \frac{1}{n} \sum_{\ell=1}^n f(\hat{X})_\ell f(\Bar{X}^i_{\hat{X}})_\ell,
\end{aligned}
\label{46}
\end{equation}
and
\begin{equation}
    \bigg( \int \mathbb{E}_{X_{\sim i}} (Y|X_i) dX_i\bigg)^2 = f^2(X) = \frac{1}{n} \sum_{\ell=1}^n f(\Bar{X})_\ell f(\hat{X})_\ell.
    \label{47}
\end{equation}
Using (\ref{46}) and (\ref{47}), we define the local sensitivity index as \cite{Saltelli10}
\begin{equation*}
    S_i = \frac{\mathrm{Var}_{X_i} \Big( \mathbb{E}_{X_{\sim i}} (Y | X_i) \Big)}{\mathrm{Var}(Y)} = \frac{ \frac{1}{n} \sum_{\ell=1}^n f(\hat{X})_\ell \Big( f(\Bar{X}^i_{\hat{X}})_\ell - f(\Bar{X})_\ell \Big)}{\mathrm{Var}(Y)}.
\end{equation*}
The global sensitivity index $S_{T_i}$ is derived in a similar way, see 
\cite{Jansen99,Saltelli10} as
\begin{equation*}
    S_{T_i} = \frac{\mathbb{E}_{X_{\sim i}} \Big( \mathrm{Var}_{X_{i}} (Y | X_{\sim i}) \Big)}{\mathrm{Var}(Y)} = \frac{ \frac{1}{2n} \sum_{\ell=1}^n \Big( f(\Bar{X})_\ell -  f(\Bar{X}^i_{\hat{X}})_\ell \Big)}{\mathrm{Var}(Y)},
\end{equation*}
where $f(\Bar{X}),f(\hat{X})$, and $ f(\Bar{X}^i_{\hat{X}})$ are the model output results using the parameter sample values from matrices $\Bar{X},\hat{X}$, and $\Bar{X}^i_{\hat{X}}$, respectively.

The independent distribution of each numerical error is determined as follows. The discretization error $\epsilon_h$ depends on the grid size $h$ used to construct the full model. We sample different grid sizes with an appropriate grid refinement ratio and construct different full models. We further calculate the discretization error associated with each model and use generated samples to construct the distribution for $\epsilon_h$.

Similar to the discretization error, the projection error $\epsilon_{\mathrm{POD}}$ and the reduced model error $\epsilon_{\mathrm{RM}}$ can be regulated by varying the reduced dimension $d$.
We have to solve the full model for some training parameter groups to construct the snapshot matrix in order to generate the reduced basis.
Before obtaining the distributions of the projection error and the reduced model error, the full model dimension and the number of training parameters are fixed.
We construct the reduced models with different reduced dimensions and obtain the corresponding errors $\epsilon_{\mathrm{POD}}$, and $\epsilon_{\mathrm{RM}}$. These error samples are then used to construct the desired distributions.

The training parameters used to construct the surrogate model govern its quality. We randomly sample the parameter groups that initiate the adaptive algorithm (parameter set $\hat{\mathcal{P}}_0$) as well as vary the cardinality of set $\hat{\mathcal{P}}$. The first few samples of the surrogate model error are obtained by randomly sampling the parameter set $\hat{\mathcal{P}}_0$ while keeping the cardinality of $\hat{\mathcal{P}}$ constant. Later on, the cardinality of $\hat{\mathcal{P}}$ is varied while $\hat{\mathcal{P}}_0$ is fixed by specifying the particular seed number.

The obtained samples are then randomly separated into two matrices $\Bar{X}$, and $\hat{X}$. Each column of the matrices $\Bar{X}$ or $\hat{X}$ represents the distribution of one numerical error, while each row describes a set of model factors, i.e., in this case, the numerical error $\{\epsilon_h, \epsilon_{\mathrm{POD}}, \epsilon_{\mathrm{RM}}, \epsilon_{samp}\}$. The total error (\ref{11}) presented in Section \ref{ErrAnalysis} is then used as a model output to perform the sensitivity analysis based on the generated sample data.
\subsection{Model Hierarchy: Stopping Criteria}
\label{ModelHierarchy}
\begin{figure}[htb]
\begin{center}
\usetikzlibrary{fit}
\usetikzlibrary{shapes,arrows,positioning,decorations.pathreplacing}
\tikzstyle{decision} = [diamond, draw, minimum height=3em, text centered]
\tikzstyle{line}=[draw]
\tikzstyle{block} = [rectangle, draw,
    text width=8cm, text centered, rounded corners, minimum height=3em]
\tikzstyle{line} = [draw, -latex']
\tikzstyle{container} = [draw, rectangle, dashed, inner sep=0.3cm]
\tikzstyle{io} = [trapezium, trapezium left angle=70, trapezium right angle=110, minimum height=3em, text centered, draw=black]
\begin{tikzpicture}[node distance = 2.5cm, auto]
\node[io, text width=3cm](Input){Parameter space $\mathcal{P}$ (of $\mathcal{N}$ scenarios), instrument, $e_{tol}$};
\node[block, below of=Input, text width=3cm](PDE){Model selection:\\PDE};
\node[block, right of=PDE, text width=3cm, xshift=2.5cm](Disc){Discretization};
\node[decision,text width=0.5cm, right of=Disc, xshift=2.5cm](TolCheck1){$\epsilon_h \leq e_{tol}^h$};
\node[block, above of=TolCheck1, text width=2cm](Refine){refine};
\node[block, below of=PDE, text width=3cm](CG){Classical greedy sampling};
\node[block, right of=CG, text width=3cm, xshift=2.5cm](POD){POD:\\RB with $\ell$ parameter groups};
\node[block, right of=POD, text width=3cm, xshift=2.5cm](RM){Reduced model};
\node[decision, below of=POD, text width=0.5cm](TolCheck2){$\epsilon_d \leq e_{tol}^d$};
\node[block, left of=TolCheck2, text width=3cm, xshift=-2.5cm](AG){Adaptive greedy sampling};
\node[block, below of=TolCheck2, text width=2cm](Stop){Stop};
\node[block, below of=AG, text width=3cm](MOR){Model order reduction};
\node[block, right of=Stop, text width=3cm, xshift=2.5cm](ScenarioCal){Solve $\mathcal{N} - \ell$ scenarios using RM};
\node[block, below of=ScenarioCal, text width=3cm](Sort){Sort $10\,000$ values};
\node[block, below of=Stop, text width=3cm](ScenarioRes){Compute $V_{fav}$, $V_{mod}$ and $V_{unfav}$};
\node[block, below of=MOR, text width=2cm](Stop2){Stop};
\path[line] (Input) -- (PDE);
\path[line] (PDE) -- (Disc);
\path[line] (Disc) -- (TolCheck1);
\path[line] (TolCheck1) --node[anchor=west] {false} (Refine);
\draw[->] (Refine) -| (Disc);
\draw[->] (TolCheck1) |- ++(0,-1.2cm) node[below] {true} -| (CG);
\path[line] (CG) -- (POD);
\path[line] (POD) -- (RM);
\draw[->] (RM) |- (TolCheck2);
\path[line] (TolCheck2) -- node[above] {false} (AG);
\path[line] (AG) -- (MOR);
\path[line] (MOR) -- (Stop);
\path[line] (TolCheck2) -- node[anchor=west] {true} (Stop);
\path[line] (Stop) -- (ScenarioCal);
\path[line] (ScenarioCal) -- (Sort);
\path[line] (Sort) -- (ScenarioRes);
\path[line] (ScenarioRes) -- (Stop2);
\end{tikzpicture}
\end{center}
\caption{Model hierarchy for model order reduction}
\label{fig:2}
\end{figure}
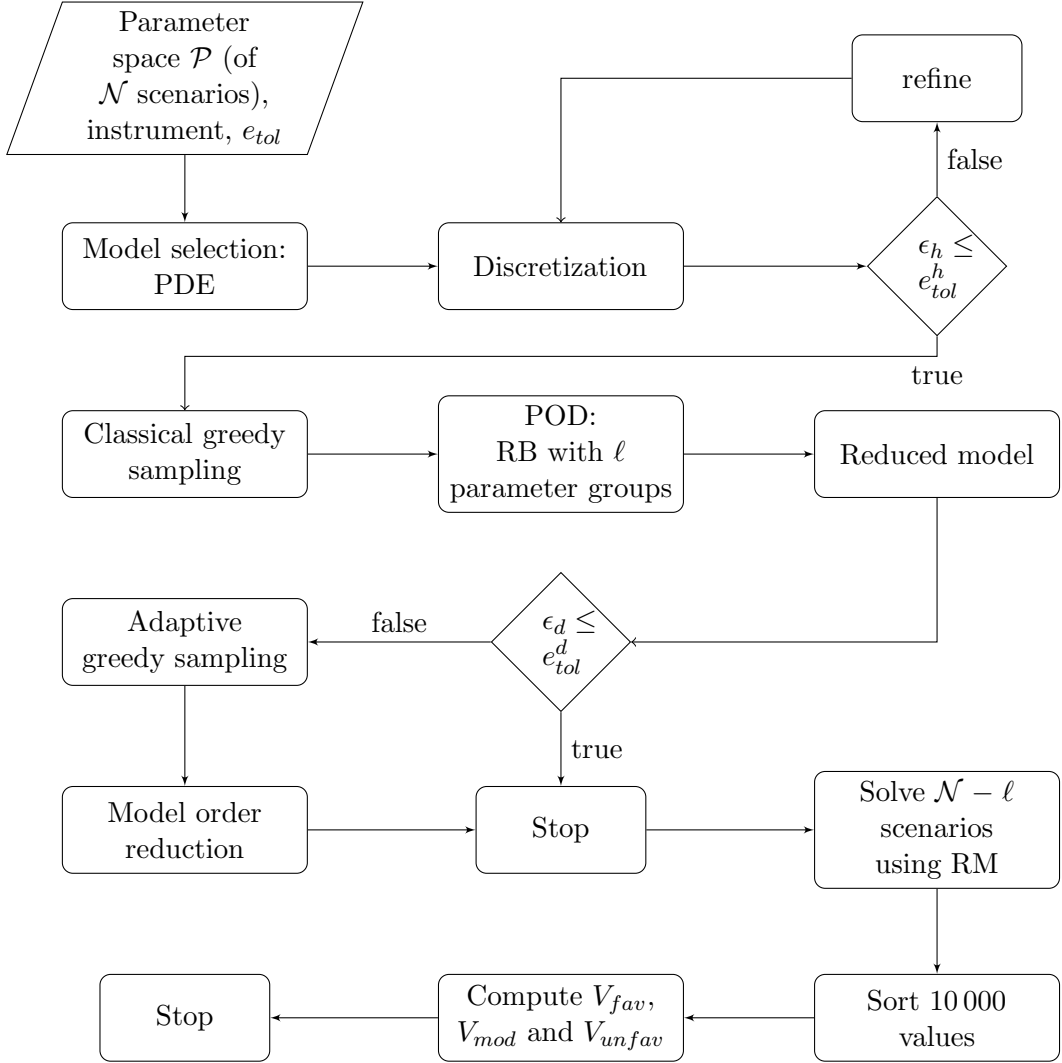

We now use the described techniques to calculate the financial risk associated with an invested product. Underlying risk factors such as interest rates have a direct influence on the invested asset. Interest rate risk arises when unanticipated developments shift the yield curve or change its shape. Thus, to calculate the financial risk, the instrument is evaluated for several thousand different simulated yield curves. One can consider these simulated yield curves as different scenarios.
The financial instruments are evaluated via the dynamics of short-rate models $\mathcal{M}$, based on the convection-diffusion-reaction partial differential equations (\ref{2}). The choice of the short-rate model depends on the underlying financial instrument.
The model hierarchy simplifies the process of obtaining a reduced model and is depicted in Fig. \ref{fig:2}.
It starts by discretizing the partial differential equation using either the finite difference method or the finite element method. The discretized model is the full model $\mathcal{M}_h$. To obtain a reduced model $\mathcal{M}_h^d$, we compute a reduced basis via the POD approach along with the classical or adaptive greedy sampling. Finally, the full model is projected onto the reduced basis to get the reduced model. The above-described steps finally lead to the desired reduced model.
However, it may happen that the user does not have enough knowledge or background about these methods to obtain the optimal reduced model. Thus, the user prefers an automatic approach that produces a reduced model whose error is below a user-defined tolerance. We describe the methodology to obtain a reduced model as a black box for which the user only requires to provide inputs and gets the desired output.
The black box method automatically selects the partial differential equation according to the underlying financial instrument, generates the discretized model, obtains the training parameters, generates the desired reduced basis, and produces a reduced model with the desired user-defined tolerance. Yield curves and a financial instrument are fed as inputs to the black box, which outputs the scenario results.
The procedure is illustrated in Fig. \ref{fig:2}.
This approach also yields a well defined stopping criteria as follows.

Consider a parameter space $\mathcal{P} \in \mathbb{R}^{\mathcal{N}\times m}$ (of scenarios $\mathcal{N}$) for which the underlying financial instrument is to be solved using an appropriate short-rate model $\mathcal{M}$. Let $\ell \in \mathcal{P} \subset \mathcal{N}$ be the number of training parameters selected using a sampling algorithm. Now, these $\ell$ parameters are used to obtain the reduced model $\mathcal{M}_h^d$ with the total error $\epsilon_T \leq \epsilon_h + \epsilon_d + \epsilon_{samp}$.
\begin{remark}{\rm 
The stopping criteria to valuate an underlying financial instrument with $\mathcal{N}$ scenarios states that when $\epsilon_T \leq e_{tol}$, where $e_{tol}$ is the user-defined tolerance, only $\ell$ scenarios need to be evaluated using the full model $\mathcal{M}_h$, while the remaining $\mathcal{N} - \ell$ scenarios can be evaluated using the reduced model $\mathcal{M}_h^d$.}
\label{Remark31}
\end{remark}
We employ the $\ell$ scenarios with the full model $\mathcal{M}_h$ to obtain the reduced basis and the remaining scenarios $\mathcal{N} - \ell$ using the reduced model $\mathcal{M}_h^d$. The results obtained for all scenarios are then sorted into favorable $V_{fav}$, moderate $V_{mod}$, and unfavorable $V_{unfav}$ performance scenarios.

\section{Numerical Example}
\label{NE}
Steepeners are the financial instruments where the cashflows depend on differences in the evolution of two different reference rates with different tenors in one currency \cite{Binder013}. We consider a callable/puttable steepener instrument whose coupons depend on the difference between two \emph{constant maturity swap} (CMS) rates, i.e., CMS10 - CMS2 \cite{Brigo06}. Table \ref{tab:2} shows the properties of a puttable steepener. A puttable steepener gives the right to the investor to claim an early redemption on certain dates (typically coupon dates) for a price fixed (typically the nominal value) in the term sheet of the bond. These early exercise dates are before the maturity of the bond. In short, a put is an option that gives the right to a buyer to sell the underlying asset at an agreed price at some time. The selling price is known as a strike price. The payoff of an option is its value at the time of its exercise. Consider $K$ as the strike price of a put option and an underlying instrument with a value $V$ and maturity $T$. The payoff of the put option is
\begin{equation*}
\begin{aligned}
    P_V &=
    \begin{cases}
    K-V, & \text{if} \hspace{0.2cm} V<K.\\
    0, & \text{if} \hspace{0.2cm} V \geq K.
    \end{cases} \\
    P_V &= \mathrm{max}(K-V,0) = (K-V)^+.
\end{aligned}
\end{equation*}
\begin{table}[htb]
\caption{Numerical Example of a floater with cap and floor.}
\label{tab:1}
\setlength{\tabcolsep}{1cm}
\begin{center}
\begin{tabular}{ll}
\hline\noalign{\smallskip}
Maturity & 10 years\\
Currency & EURO \\
Coupon frequency & Annually\\
Coupons: & \\
Year 1 to Year 3 & 4\% fix \\
Year 4 to Year 10 & CMS10 - CMS2 \\
Cap rate, $C_R$ & 3.0 \% p.a. \\
Floor rate, $F_R$ & 0.0 \% p.a. \\
\noalign{\smallskip}\hline
\end{tabular}
\end{center}
\end{table}
Table \ref{tab:2} shows that the interest rate is capped at 3.0\% p.a. and floored at 0.0\% p.a. The cap is nothing but an upper limit on the interest rate for a floating interest rate instrument, while the floor sets a lower limit on the interest rate. The coupon rate from the fourth year until maturity is
\begin{equation}
\begin{aligned}
        \mathrm{coupon \; rate} =& \mathrm{min}(\mathrm{coupon \; cap}, \mathrm{max}(\mathrm{coupon \; floor}, \mathrm{CMS10} - \mathrm{CMS2})),\\
        =& \mathrm{min}(3.0\%, \mathrm{max}(0.0\%,\mathrm{coupon \; floor}, \mathrm{CMS10} - \mathrm{CMS2}))
\end{aligned}
\label{48}
\end{equation}
We solve the puttable steepener example using the two-factor Hull-White model. The partial differential equation is discretized using the finite element method to obtain the full model, which is further approximated to achieve the reduced model. The model parameters are $\sigma_1 = 0.0035$, $\sigma_2 = 0.008$, $\gamma = 0.65$, $\alpha = 0.75$, $b = 0.04$.
We have computed the model parameter $\theta_\ell$ as explained in \cite{MathConsult09}. The yield curve simulation is the first step to compute the model parameters. Based on the yield curve simulation procedure described in \cite{binder2020}, we have performed the bootstrapping process for the recommended holding period and the intermediate holding period, i.e., for ten years and five years, respectively. The collected historical data has $21$ tenor points and $1306$ observation periods as follows (D: Day, M: Month, Y: Year):
\begin{equation*}
\begin{aligned}
    m &=: \{1D, 3M, 6M, 1Y, 2Y, 3Y, \cdots, 10Y, 12Y, 15Y, 20Y, 25Y, 30Y, 40Y, 50Y\}\\
    n &=: \{\text{1306 daily interest rates at each tenor point}\}
\end{aligned}
\end{equation*}
\begin{figure}[htb]
  \centering
  \includegraphics[width=1\columnwidth]{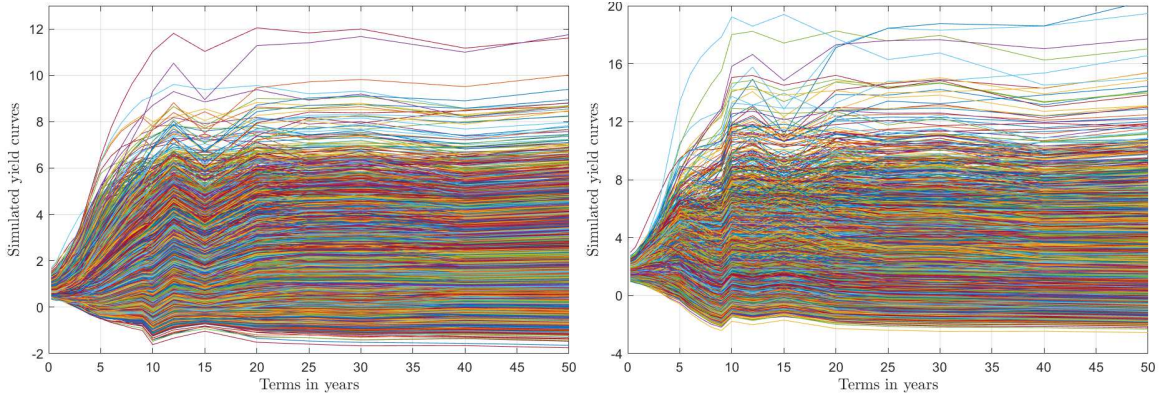}
  \caption{$10\,000$ simulated yield curves obtained by bootstrapping for five years and ten years in the future.}
  \label{fig:3}
\end{figure}
Figure \ref{fig:3} shows the plots of simulated yield curves obtained by bootstrapping for five years and ten years in the future. The calibration generates a parameter space $\mathcal{P}$ of order $10\,000 \times m$ for which we now solve the two-factor Hull-White model.

For calculations, we have used a structured, two-dimensional triangular mesh grid that generates a full model of order $M$.
We have used fully implicit time-stepping so that discretizations in time and space can be chosen independently. The resulting system of linear equations is
\begin{equation}
    A(\rho_\ell(t))V^{n-1} = B(\rho_\ell(t))V^n,
    \label{49}
\end{equation}
To obtain a solution, we solve this system starting at $t=T$ with an appropriate terminal condition $V(T)$ backward in time at each time step $n$. It is important to note that we have selected the time step $\Delta t= 20$ days using Algorithm \ref{Algo1}, and ensuring that all key dates are achievable, i.e., coupon dates, put dates, or valuation dates. The tolerance $e_{tol}^t$ of $10^{-3}$ is used to obtain this best suitable time step $\Delta t$.
We determine the value of the steepener by solving the linear system given in (\ref{49}) as follows. If the coupon date is reached, then we update the value of the steepener by adding the necessary coupon given by (\ref{48}) \cite{MathConsult09}
\begin{equation*}
    V^{n-1} = V^{n-1} + C_F^{n-1},
\end{equation*}
where $V^{n-1}$ is the value calculated by solving the full model, and $C_F$ is either an already known cashflow or a cashflow generated by the coupon if the coupon date is achieved at $t^{n-1}$ or zero. Similarly, if the put date is reached, then we have to update the value at that time point by the following max function
\begin{equation*}
    V^{n-1} = \mathrm{max}(V^{n-1},\mathrm{Put \; Price}),
\end{equation*}
where the put price is the strike price defined by the put option contract. In this paper, we have considered $\mathrm{Put \; Price}=1$. We have considered put dates annually starting one year after the issuing of the bond until one year before maturity.

\begin{figure}[htb]
  \centering
  \includegraphics[width=0.4\columnwidth]{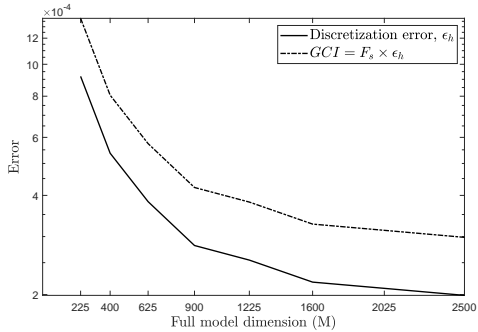}
  \caption{Discretization error estimator and  grid convergence index (\ref{26}) vs  full model size.}
  \label{fig:4}
\end{figure}

The quality of the full model is assured by calculating the discretization error associated with the two-factor Hull-White model, as explained in Subsection \ref{DiscPDE}.
Figure \ref{fig:4} shows the discretization error along with the grid convergence index values. We can see that the error decreases with increasing order of the full model.
A discretization error computation requires three different solutions of the full model for three different grid sizes. The error estimator is acceptable if these solutions are in the asymptotic range of convergence, i.e., the observed order of accuracy matches the formal order of accuracy, and the grid convergence index is low. The formal order of accuracy depends on the order of the polynomial used to define the elements in the finite element method \cite{Hutton04}. The $p$th order polynomial gives the $p+1$ order of convergence \cite{Babuska94}. We have constructed our full model using the 1st order linear triangular elements. Thus, the formal order of accuracy $p_f=2$. The observed order of accuracy $\hat{p}$ is in the range of 1.8959 to 2.5616 with an average of 2.072. For the full model with $M=1600$, $\hat{p} = 1.9213$, which matches the formal order of accuracy within $10\%$. Thus, we can say that the solutions obtained using the full model are in the asymptotic range. The safety factor used to calculate the grid convergence index (\ref{26}) is then $1.5$. The full model size is selected using Algorithm \ref{Algo5}, which is terminated when the discretization error is less than the user-defined tolerance $e_{tol}^h = 5\times 10^{-4}$. We have designed a full model of order $M= 1600$ as the discretization error $\epsilon_h = 2.182\times 10^{-4}$ falls below $5\times 10^{-4}$. Also, one can notice that the grid convergence index is less than $0.05\%$.
\subsection{Model Order Reduction}
\label{MORExample}
For the model order reduction approach, we have obtained the training parameter set to compute the optimal reduced basis by using either the classical greedy sampling or the adaptive greedy sampling. In the following, we present the results and error analysis of these algorithms. For the classical greedy sampling, we have used $c=40$ randomly selected parameter groups to train the algorithm. The allowed maximum iterations $I_{max} = 10$, i.e., at most $10$ parameter groups are selected which maximize the residual error. The algorithm then updates the snapshot matrix and generates a new reduced basis.

\begin{figure}[htb]
  \centering
  \includegraphics[width=0.4\columnwidth]{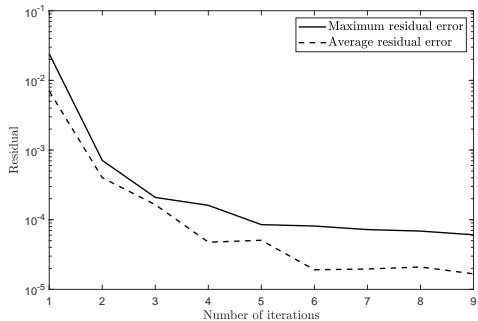}
  \caption{Evolution of maximum and average residuals with each iteration of classical greedy algorithm.}
  \label{fig:5}
\end{figure}
Figure \ref{fig:5} shows the evolution of maximum and average residuals with each iteration of the classical greedy sampling algorithm. It is observed that the residual error decreases with each proceeding iteration. We can truncate the greedy iteration after 4th or 5th iteration as the error has dropped below $10^{-4}$.

\begin{figure}[htb]
  \centering
  \includegraphics[width=0.6\columnwidth]{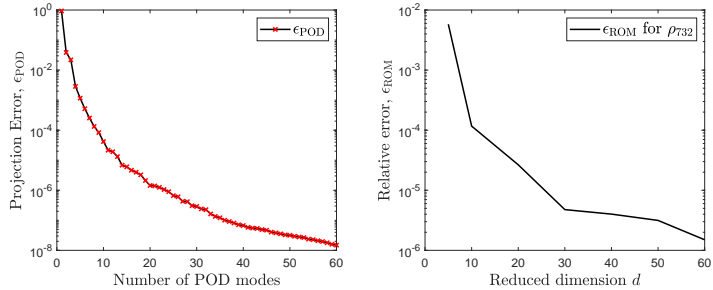}
  \caption{Projection error associated with POD (left), and a plot of relative error between  full and reduced model obtained using the classical greedy sampling approach (right).}
  \label{fig:6}
\end{figure}

The dimension of the reduced basis $d$ is chosen using Algorithm \ref{Algo5}. The projection error is plotted in Fig. \ref{fig:6}. We observe that the POD projection error decreases monotonically with increasing $d$, which is determined in terms of POD eigenvalues.
The graph of the projection error shows that we have succeeded in determining an optimal reduced basis, which is later used to generate the reduced model.
Figure \ref{fig:6} also shows the reduced model error $\epsilon_{\mathrm{RM}}$ plotted against the reduced dimension $d$ for the parameter group $\rho_{732}$. It is observed that the error $\epsilon_{\mathrm{RM}}$ decreases as the dimension $d$ increases. This shows that the reduced model is an excellent approximation of the full model. We have designed a reduced model of dimension $d=10$ as the model order reduction error $\epsilon_{d} = \epsilon_{\mathrm{POD}} + \epsilon_{\mathrm{RM}} = 4.18\times 10^{-5} + 1.17 \times 10^{-4} = 1.588 \times 10^{-4}$ is less than $e_{tol}^d = 5\times 10^{-4}$. Finally, we used this reduced model to solve the remaining scenarios, as explained in Section \ref{ModelHierarchy}.
\begin{figure}[htb]
  \centering
  \includegraphics[width=0.4\columnwidth]{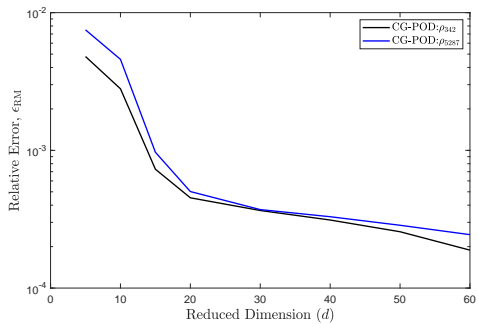}
  \caption{Drawback of classical greedy (CG) sampling approach.}
  \label{fig:7}
\end{figure}
In classical greedy sampling, the algorithm is trained on the randomly generated parameter set $\hat{\mathcal{P}}$ of cardinality $40$. We noticed that the random sampling often neglects some of the vital parameter groups within the parameter space.
The reduced model designed using the classical greedy sampling approach shows this drawback as illustrated in Fig. \ref{fig:7}.
The figure shows the reduced model error for two different parameter groups ($\rho_{342},\rho_{5287}$). In both cases, the relative error does not drop below the user-defined tolerance $e_{tol}^d =5 \times 10^{-4}$, and Algorithm \ref{Algo5} does not converge.
This prompted us to implement an adaptive greedy approach that relies on the surrogate model to locate the training parameters.
\begin{figure}[htb]
  \centering
  \includegraphics[width=0.6\columnwidth]{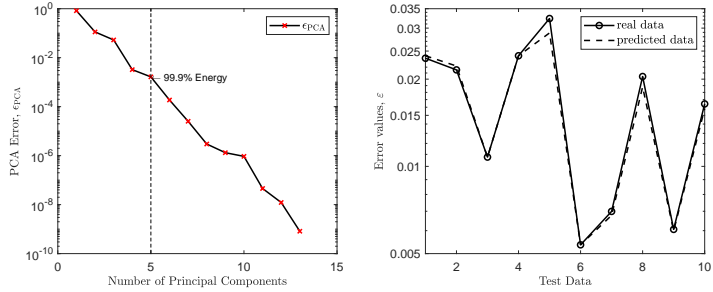}
  \caption{PCA error (left) and predicted values obtained using the surrogate model for an error estimator vs actual values (right).}
  \label{fig:8}
\end{figure}
\begin{table}[htb]
\caption{Mean squared error of prediction}
\label{tab:2}
\setlength{\tabcolsep}{1cm}
\begin{center}
\begin{tabular}{ll}
\hline\noalign{\smallskip}
Fold No. & Squared error\\
\hline\noalign{\smallskip}
$L_K = 1$ & $8.7790 \times 10^{-5}$ \\
$L_K = 2$ & $7.6331 \times 10^{-5}$ \\
$L_K = 3$ & $7.9755 \times 10^{-5}$ \\
$L_K = 4$ & $8.8776 \times 10^{-5}$ \\
\hline\noalign{\smallskip}
$\mathrm{MSEP}_{\mathrm{K-CV}}$ & $8.3163 \times 10^{-5}$ \\
\hline\noalign{\smallskip}
$\mathrm{MSEP}_{\mathrm{K-CVadj}}$ & $7.2815 \times 10^{-5}$ \\
\noalign{\smallskip}\hline
$\epsilon_{samp} = \mathrm{MSEP}_{\mathrm{K-CVadj}} + \epsilon_{\mathrm{PCA}}$ & $2.3750\times 10^{-4}$ \\
\noalign{\smallskip}\hline
\end{tabular}
\end{center}
\end{table}
The used surrogate model  is based on the principal component regression technique \cite{binder2020}.
The adaptive greedy iteration is initiated by computing the error estimator for randomly selected parameter set $\hat{\mathcal{P}}_0$ of cardinality $c_0 = 10$. At each surrogate model iteration, we have chosen $c_k = 10$ parameter groups. This process repeat itself until we get the parameter set $\hat{\mathcal{P}}$ of cardinality 40.
Figure \ref{fig:8} shows monotonically decreasing singular values of the data matrix $\hat{\mathcal{P}}_k$ composed of $c_k$ parameter groups at $k$th iteration.
We noticed that the first three principal components comprise almost 99\% of energy, and within the first five components, energy level round up to 99.99\%. Thus, the surrogate model has only five principal components.
Figure \ref{fig:8} demonstrates that the predicted values obtained via the surrogate model match the actual values within an acceptable range.

The quality of the designed surrogate model is further tested using a mean squared error of prediction. We have used a K-fold cross-validation technique to determine the MSEP as described in Subsection \ref{SampError}. The learning data matrix composed of $k\times c_k$ parameter groups and error estimators at the $k$th iteration is divided into $K = 4$ folds. We obtain the mean square error for each fold and sum it to calculate the MSEP of the cross-validation. Table \ref{tab:2} shows the calculations and MSEP for the designed surrogate model.
The sampling error $\epsilon_{samp}$ is the sum of the adjusted K-fold cross-validation error and the projection error, i.e., $\epsilon_{samp} = \mathrm{MSEP}_{\mathrm{K-CVadj}} + \epsilon_{\mathrm{PCA}} = 2.3750 \times 10^{-4}$.
As $\epsilon_{samp} < e_{tol}^{samp} = 5\times 10^{-4}$, we can say that the designed surrogate model is satisfactory.
\begin{figure}[!htb]
  \centering
  \includegraphics[width=0.4\columnwidth]{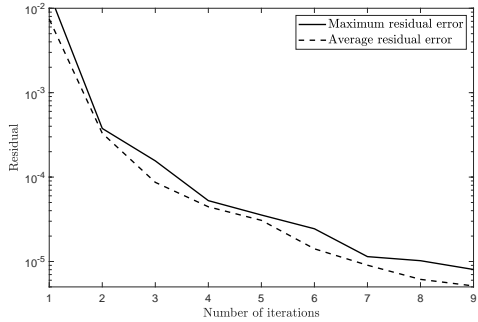}
  \caption{Evolution of maximum and average residuals in the adaptive greedy algorithm.}
  \label{fig:9}
\end{figure}
We have used this surrogate model to construct the parameter set $\hat{\mathcal{P}}$. Finally, the optimal parameter group $\rho_I$ is selected from the parameter set $\hat{\mathcal{P}}$ that maximizes the residual error. The reduced basis is obtained by computing the truncated SVD of this snapshot matrix, as shown in Algorithm \ref{Algo5}.
Figure \ref{fig:9} shows the evolution of the maximum and average residual error with each adaptive greedy iteration. It illustrates that the residual error decreases with increasing iteration and falls below $10^{-4}$ after the fourth iteration, which is enough to truncate the greedy iteration. The decrease of the residual error with each incrementing iteration also shows that the algorithm succeeded in efficiently locating the optimal parameter group.
\begin{figure}[!htb]
  \centering
  \includegraphics[width=0.8\columnwidth]{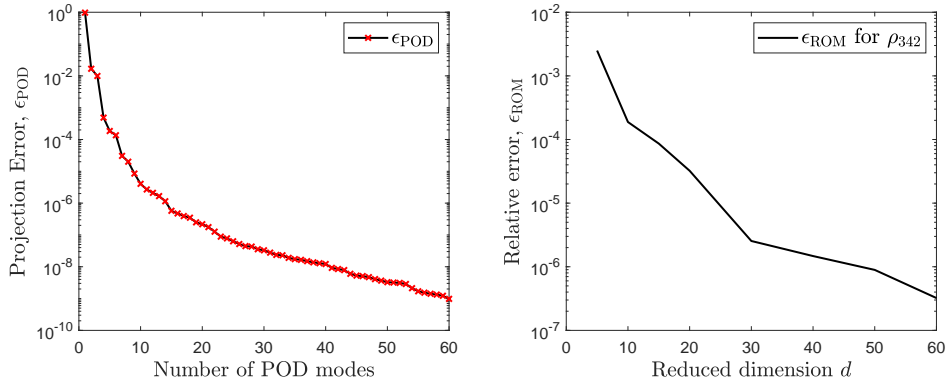}
  \caption{Projection error associated with POD (left), and plot of relative error between full and reduced model obtained using adaptive greedy sampling (right).}
  \label{fig:10}
\end{figure}
The reduced basis dimension $d$ is selected such that for the generated reduced model, $\epsilon_{\mathrm{POD}}$ and $\epsilon_{\mathrm{RM}}$ are below the user defined tolerance. Figure \ref{fig:10} shows the monotonically decreasing projection error $\epsilon_{\mathrm{POD}}$ associated with the POD. The reduced model error $\epsilon_{\mathrm{RM}}$ plotted in Fig. \ref{fig:10} shows that the error decreases with increasing the reduced dimension $d$. We have constructed the final reduced model using $d=10$, as the model order reduction error is less than the user-defined tolerance $e_{tol}^d$, i.e., $\epsilon_{d} = \epsilon_{\mathrm{POD}} + \epsilon_{\mathrm{RM}} = 4.042 \times 10^{-6} + 1.882 \times 10^{-4} = 1.922\times 10^{-4} < 5\times 10^{-4}$. This means that the reduced model is an excellent approximation of the full model.
\begin{figure}[!htb]
  \centering
  \includegraphics[width=0.4\columnwidth]{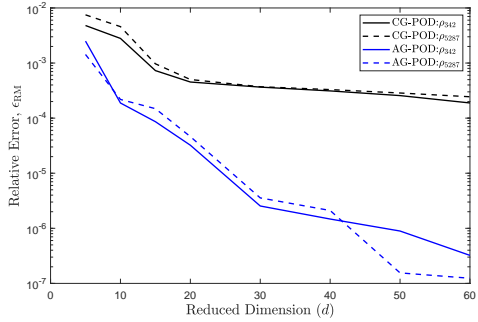}
  \caption{Comparison of classical greedy (CG) approach and adaptive greedy (AG) approach.}
  \label{fig:11}
\end{figure}
We have performed $10$ full model evaluations during the adaptive greedy sampling, and the remaining scenarios are solved using this reduced model, as explained in Subsection \ref{ModelHierarchy}.
We have implemented the adaptive greedy sampling approach to overcome the drawbacks of the classical greedy sampling method, as evident from Fig. \ref{fig:11}. The relative error between the full model and the reduced model obtained using the adaptive greedy sampling approach is smaller than $10^{-4}$ with the reduced dimension $d=10$ for both parameter groups $\rho_{342}, \rho_{5287}$. The error $\epsilon_{\mathrm{RM}}$ keeps on decreasing with increasing $d$, unlike in the classical greedy sampling approach.

\begin{figure}[!htb]
  \centering
  \includegraphics[width=0.6\columnwidth]{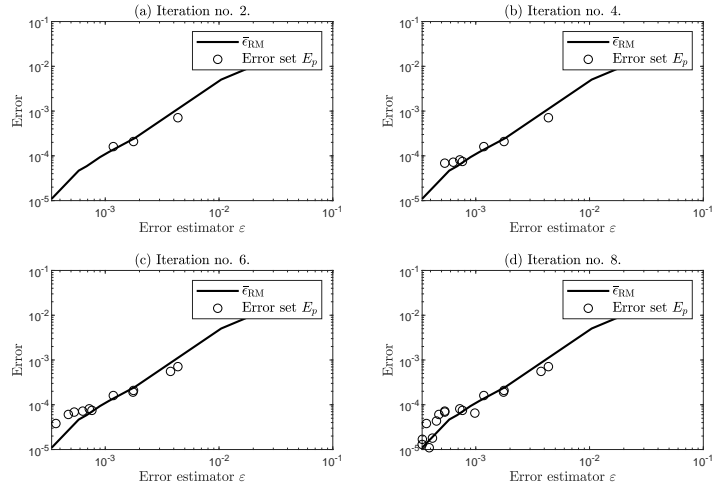}
  \caption{Error model $\Bar{\epsilon}_{\mathrm{RM}}$ based on the available error set $E_p$ for $4$ different greedy iterations.}
  \label{fig:12}
\end{figure}

The convergence of the classical greedy sampling approach is monitored using the residual error. However, we have used an error model $\Bar{\epsilon}_{\mathrm{RM}}$ that uses the reduced model errors and the residual errors to monitor the convergence of the adaptive greedy sampling. Figure \ref{fig:12} illustrates the designed error model based on the available error set $E_p$ for $4$ different greedy iterations. The error plot demonstrates a strong correlation between the relative error and the residual error. The results indicate that a consideration of the linear error model is adequate to capture the overall behavior of the exact error as a function of the residual error.
\subsubsection{Random Sampling}
\label{RandSamp}
\begin{figure}[htb]
  \centering
  \includegraphics[width=0.4\columnwidth]{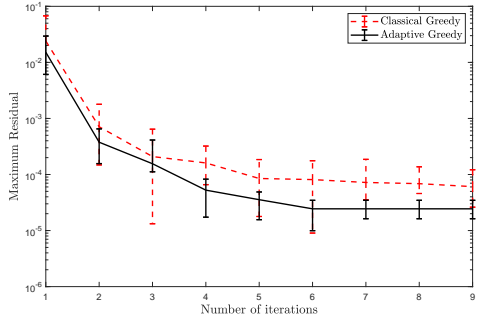}
  \caption{Effect of random sampling on classical and adaptive greedy sampling.}
  \label{fig:13}
\end{figure}
We have used random sampling up to a certain extent in both classical and adaptive greedy approaches.
Random sampling is used to construct the parameter set $\hat{\mathcal{P}}$ in the classical greedy approach, while in the case of the adaptive greedy approach, it is used to initiate the algorithm.  In both cases, the random sampling somewhat affects the results. Figure \ref{fig:13} shows the maximum residual error plot with each proceeding iteration of the greedy algorithms with error bars.
These results are obtained by performing classical and adaptive procedures for $30$ different sets of random numbers.
It is observed that the adaptive greedy approach is less sensitive to random sampling than its classical counterpart. On average, the adaptive greedy approach yields a better reduced model as its residual errors are smaller than the classical approach.
\subsubsection{Sensitivity Analysis}
\label{SAExample}
We have obtained $30$ different data points for each numerical error to generate respective distributions by varying certain key parameters.
We have sampled $30$ different grid sizes to compute the discretization error distribution, where the full model dimension is being varied from $100$ to $2500$. To compute the sampling error distribution, we have randomly generated $15$ different $\hat{\mathcal{P}}_0$ by keeping the cardinality of $n(\hat{\mathcal{P}}) = 50$. Additionally, $15$ data points for the sampling error are calculated by varying $n(\hat{\mathcal{P}})$ between $30$ and $100$, while keeping $\hat{\mathcal{P}}_0$ fixed using a particular seed number. The distributions for the projection error and the reduced model error are obtained by varying the reduced dimension between $5$ and $60$.

The local indices show that the final output is most sensitive to the discretization error followed by the sampling error, while it is least sensitive to the reduced model error and the projection error, respectively. We can mitigate the discretization error by increasing the full model dimension. Nevertheless, sensitivity indices prevail that the change in grid size considerably affects the final output. We overcome this problem by increasing the full model dimension until the further change in grid sizes does affect the solutions significantly. We have constructed the full model of order $M=1600$, where the discretization error is of order $10^{-4}$.
\begin{table}[htb]
\caption{Sensitivity analysis: local and global sensitivity indices.}
\label{tab:3}
\begin{center}
\begin{tabular}{lllll}
\hline\noalign{\smallskip}
Numerical Error & Local Index & Rank & Global Index & Rank \\
\noalign{\smallskip}\hline\noalign{\smallskip}
Dicretization error $\epsilon_h$ & 0.3351 & 1 & 0.4172 & 2 \\
Reduced model error $\epsilon_{\mathrm{RM}}$ & 0.1450 & 3 & 0.0917 & 3 \\
Projection error $\epsilon_{\mathrm{POD}}$ & 0.0167 & 4 & 0.0081 & 4 \\
Sampling error $\epsilon_{samp}$ & 0.2247 & 2 & 0.5397 & 1 \\
\noalign{\smallskip}\hline
\end{tabular}
\end{center}
\end{table}

The global sensitivity analysis shows that the sampling error is the most sensitive numerical error. The major input factors contributing to the sampling error are the randomly sampled parameter set $\hat{\mathcal{P}}_0$ and the cardinality of set $\hat{\mathcal{P}}$. As it is not feasible to search the entire parameter space as the parameter set $\hat{\mathcal{P}}$, the surrogate model only locates a certain number of parameter groups to construct $\hat{\mathcal{P}}$ and neglects the rest. This omission of  parameter groups contributes to the sampling error together with the random sampling of $\hat{\mathcal{P}}_0$ used to initiate the surrogate model loop, as we have already demonstrated in Figure 13. This random sampling contributes to the uncertainty and makes the sampling error more sensitive.
To overcome this problem, we suggest to increase the cardinality of $\hat{\mathcal{P}}$ used to locate the most optimal parameter group if necessary. We have kept the cardinality of $\hat{\mathcal{P}}$ as $40$, as the sampling error is of order $10^{-4}$.

Although increasing the cardinality of the set will increase the computational burden, it generates the best suitable reduced basis that allows constructing a fairly low dimensional reduced model. This compensates for the extra computational burden necessary during the adaptive sampling.
We noticed that the projection error is the least sensitive to the final output considering both local (individual contribution) and global (interactions) indices. This tells us that the singular value decay of the snapshot matrix is efficient to obtain the desired reduced basis.
\subsubsection{Computational Cost}
\label{CompCost}
The three major error contributors that lead to the total error are the discretization error, the model order reduction error, and the sampling error.
\begin{figure}[htb]
\begin{center}
\usetikzlibrary{shapes,arrows,positioning,decorations.pathreplacing}
\tikzstyle{block} = [rectangle, draw,
    text width=8cm, text centered, rounded corners, minimum height=3em]
\tikzstyle{line} = [draw, -latex']

\begin{tikzpicture}[node distance = 1.5cm, auto]
\node [block] (MainModel) {Partial differential equation, $\mathcal{M}$
};
\node[draw=none,fill=none, below of=MainModel] (DisError){ Discretization error, $\epsilon_h = 2.182\times 10^{-4}$};
\node[block, below of=DisError](DiscModel){Discretized model, $\mathcal{M}_h$
};
\node[draw=none,fill=none, below of=DiscModel, text width=8cm, yshift=-0.5cm](MORError){ Reduced model error, $\epsilon_{\mathrm{RM}} = 1.882\times 10^{-4}$\\
Projection error, $\epsilon_{\mathrm{POD}} = 4.042\times 10^{-6}$\\
Sampling error, $\epsilon_{\mathrm{samp.}}= 2.3750\times 10^{-4}$
};
\node[block, below of=MORError, yshift=-0.5cm](RM){Reduced model, $\mathcal{M}_h^d$
};
\path [line] (MainModel) -- (DisError);
\path [line] (DisError) -- (DiscModel);
\path [line] (DiscModel) -- (MORError);
\path [line] (MORError) -- (RM);
\draw[decorate,decoration={brace,amplitude=20pt}]
(MainModel.north east) -- (RM.south east) node [midway,xshift=5cm,left,text width=4cm] {Total error,\\$\epsilon_T \approx 6.479 \times 10^{-4}$};
\end{tikzpicture}
\end{center}
\caption{Model hierarchy showing errors arising in the analysis of the mathematical model.}
\label{fig:14}
\end{figure}
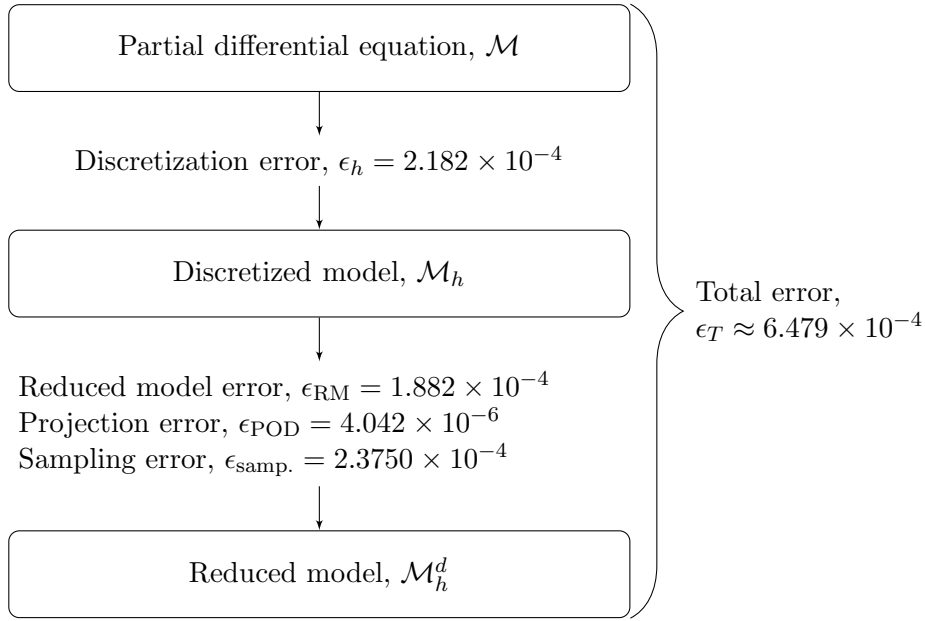
Figure \ref{fig:14} presents the model hierarchy along with the error associated with each stage.
The total error is a sum of all significant errors and is given by (\ref{11})
\begin{equation*}
\begin{aligned}
        \epsilon_T &\approx \epsilon_h + \epsilon_{\mathrm{RM}} + \epsilon_{\mathrm{POD}} + \epsilon_{\mathrm{samp.}},\\
        \epsilon_T &= 2.182\times 10^{-4} + 1.882\times 10^{-4} + 4.042 \times 10^{-6} + 2.3750\times 10^{-4} = 6.479 \times 10^{-4}.
\end{aligned}
\end{equation*}
As $\epsilon_T < 10^{-3}$ fulfills the stopping criteria prescribed in Remark \ref{Remark31}. Thus, we can say that $\ell = 10$ full model valuations are enough to obtain a fairly accurate reduced model, and the remaining scenarios (parameter groups) can be solved using the generated reduced model.

Calculation of the computation time is a crucial aspect enlightening the importance of the model order reduction over the full model approach. In the classical greedy sampling approach, the algorithm solves $c$ reduced models and one full model at each greedy iteration. The algorithm then updates the snapshot matrix using the full model solution and computes its truncated SVD. Let $t_{\mathrm{RM}}$ be the time required to solve one reduced model, $t_{\mathrm{FM}}$ be the computational time required for one full model, and let $t_{\mathrm{SVD}}$ be the time required to obtain a truncated SVD of the snapshot matrix. The total computational time $T_{Q}^{\mathrm{CG}}$ required to obtain the reduced basis after $i$ iterations of the classical greedy sampling approach can then be given as
\begin{equation*}
    T_{Q}^{\mathrm{CG}} \approx \bigg [ c \times t_{\mathrm{RM}} +  (t_{\mathrm{FM}} + t_{\mathrm{SVD}}) \bigg] \times i.
\end{equation*}
The adaptive greedy approach at each iteration conducts the following steps to obtain a reduced basis. It is initiated by solving $c_0$ reduced models to obtain the residual errors $\{\varepsilon_{j}\}_{j=1}^{c_0}$. These error estimator values are used to construct a surrogate model that locates $c_k$ parameter groups. The process is repeated until the total cardinality of $\hat{\mathcal{P}}$ reaches $c$, i.e., for $k$ iterations. The algorithm then solves the full model for the optimal parameter group, updates the snapshot matrix, and computes a truncated singular value decomposition of the updated snapshot matrix. Finally, an error model $\Bar{\epsilon}_{\mathrm{RM}}$ is built to monitor the convergence of the greedy algorithm. Thus, the total computational time can be given as
\begin{equation*}
    T_Q^{\mathrm{AG}} \approx \bigg[ c_0 \times t_{\mathrm{RM}} + k (c_k \times t_{RM} + t_{\mathrm{SM}} + t_{\mathrm{SM}}^{ev}) + t_{\mathrm{FM}} + t_{\mathrm{SVD}} + 2t_{\mathrm{RM}}^{aft,bef} + t_{\mathrm{EM}} \bigg ] \times i,
\end{equation*}
where $t_{\mathrm{SM}}$ and $t_{\mathrm{SM}}^{ev}$ denote the computational times required to build and valuate a surrogate model for the entire parameter space, respectively where $t_{\mathrm{EM}}$ is the time required to build an error model. The term $2t_{\mathrm{RM}}^{aft,bef}$ shows the computational time needed to solve the reduced model after and before updating the reduced basis.
\begin{table}[htb]
\caption{Computational time comparison: SVD vs Randomized SVD.}
\label{tab:4}
\begin{center}
\begin{tabular}{cccc}
\hline\noalign{\smallskip}
 & SVD & RandSVD & Speedup factor\\
\noalign{\smallskip}\hline\noalign{\smallskip}
Computational time $t_{\mathrm{SVD}}$ & 8.756 s & 0.7102 s & $\approx$ 12.5 \\
\noalign{\smallskip}\hline
\end{tabular}
\end{center}
\end{table}
We have used a truncated SVD based on randomized algorithms to obtain a reduced basis. Table \ref{tab:4} shows the computation time comparison between the randomized truncated and full SVD. For our problem, the randomized SVD is at least $12$ times faster than the basic SVD.
\begin{table}[htb]
\caption{Computing time/ reduction time ($T_{Q}$) to generate a reduced basis.}
\label{tab:5}
\begin{center}
\begin{tabular}{llll}
\hline\noalign{\smallskip}
Algorithm & Cardinality $|\hat{\mathcal{P}}|$ & $I_{max}$ & Computing time \\
\noalign{\smallskip}\hline\noalign{\smallskip}
Classical greedy sampling & 20 & 10 & 106.81 s \\
Classical greedy sampling & 30 & 10 & 156.24 s \\
Classical greedy sampling & 40 & 10 & 278.46 s \\
Adaptive greedy sampling & 40 & 10 & 387.63 s \\
\noalign{\smallskip}\hline
\end{tabular}
\end{center}
\end{table}
\begin{table}[htb]
\caption{Evaluation time.}
\label{tab:6}
\begin{center}
\begin{tabular}{lllll}
\hline\noalign{\smallskip}
Algorithm & Model & \makecell{Eva. time \\single $\rho_s$} & \makecell{Total Eva.\\ time ($T_{\mathrm{eva}}$)} & \makecell{Total time \\ $T_Q + T_{\mathrm{eva}}$} \\
\noalign{\smallskip}\hline\noalign{\smallskip}
 & FM, $M=1600$ & 1.9136 s & 19316.5 s & 19316.5 s \\
\noalign{\smallskip}\hline\noalign{\smallskip}
Classical greedy sampling & RM, $d = 5$ & 0.198 s & 1975.31 s & 2253.7 s \\
Classical greedy sampling & RM, $d = 10$ & 0.268 s & 2681.45 s &  2959.9 s  \\
\noalign{\smallskip}\hline\noalign{\smallskip}
Adaptive greedy sampling & RM, $d= 5$ & 0.214 s & 2143.47 s & 2531.1 s\\
Adaptive greedy sampling & RM, $d= 10$ & 0.263 s & 2671.30 s & 3058.9 s\\
\noalign{\smallskip}\hline
\end{tabular}
\end{center}
\end{table}
We noticed that the evaluation of the reduced model is at least $8-10$ times faster than the full model. However, there is a slight increase in the total time due to the additional reduction time. Table \ref{tab:6} also infers that the evaluation time increases with an increasing reduced dimension $d$. The reduced model obtained based on the classical greedy sampling approach is at least $7-9$ times faster than the full model. The computational time for the adaptive greedy based reduced model is a bit higher due to the time invested in computing surrogate and error models. The reduced model obtained using the adaptive greedy sampling approach is $6-8$ times faster than the full model. The computation time presented in the table considers that the greedy algorithms run for the maximum number of iterations $I_{max}$. However, in practice, we can truncate the algorithms after  the $4$th or $5$th iteration and achieve even more speedup.
\subsection{Steepener Results}
\label{SteepenerRes}
We solve the steepener instrument with the reduced model obtained using the adaptive greedy approach for $10\,000$ different parameter groups. To design a key information document, we need the values of the steepener at different spot rates. The spot rate is nothing but the yield rate at the first tenor point from the simulated yield curve. For $10\,000$ simulated yield curves, we procure $10\,000$ different spot rates and the corresponding values for the instrument.
\begin{table}[htb]
\caption{Results for the puttable steepener.}
\label{tab:7}
\setlength{\tabcolsep}{0.4cm}
\begin{center}
\begin{tabular}{lllll}
\hline\noalign{\smallskip}
Performance Scenario & 5 years &  & 10 years & \\
\noalign{\smallskip}\hline
Model & RM & UnRisk & RM & UnRisk\\
\noalign{\smallskip}\hline
Favorable (90th percentile) & 0.983 & 0.972 & 1.001 & 0.994 \\
Moderate (50th percentile) & 0.931 & 0.926 & 0.940 & 0.934 \\
Unfavorable (10th percentile) & 0.907 & 0.904 & 0.912 & 0.919 \\
\noalign{\smallskip}\hline
\end{tabular}
\end{center}
\end{table}
The regulations also demand to find a VaR equivalent volatility (VEV) to determine a market risk indicator at the recommended holding period.
The VaR shall be the instrument's price at a confidence level of $97.5\%$ at the end of the recommended holding period discounted to the present date using the risk-free discount factor $\mathrm{(DF)}$. Thus, the VaR in price space is $\mathrm{VaR}\times \mathrm{DF}$, and the VEV is
\begin{equation*}
    \mathrm{VEV} = \frac{ \sqrt{(3.842 - 2\times \mathrm{ln}(\mathrm{VaR}_{\mathrm{price\; space}}))} - 1.96 }{\sqrt{T}}.
\end{equation*}
\begin{figure}[htb]
  \centering
  \includegraphics[width=1\columnwidth]{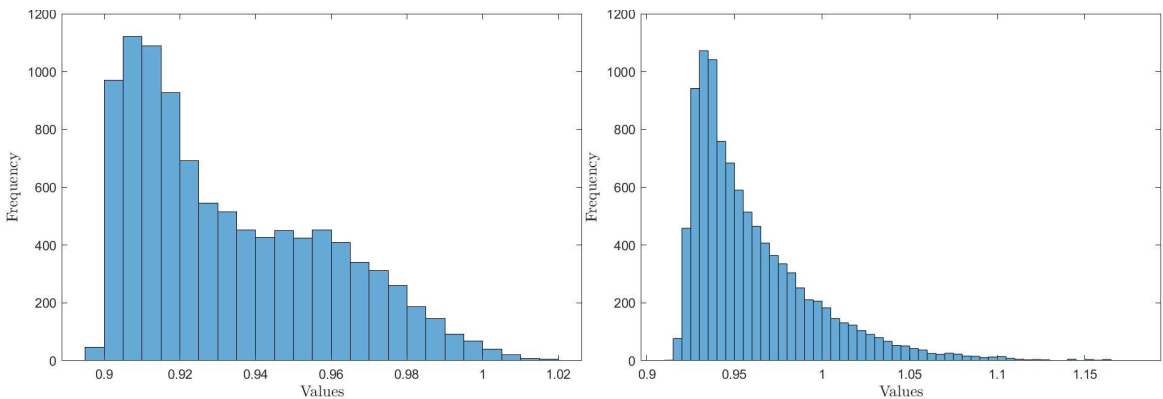}
  \caption{Distribution of $10\,000$ results after five years (left) and ten years (right).}
  \label{fig:15}
\end{figure}
For the steepener instrument, the calculated VEV is $1.3285 \%$ which gives the market risk indicator as $2$. See \cite{binder2020} for a detailed discussion on the VEV and market risk indicators. It is required to include three different performance scenarios in the KID: (i) favorable scenario, (ii) moderate scenario, (iii) unfavorable scenario, which are the values at 90th percentile, 50th percentile, 10th percentile of $10\,000$ values, respectively. The steepener performance scenarios obtained using the reduced model and the commercially available software UnRisk for the comparison are presented in the table \ref{tab:7} after five and ten years.
Figure \ref{fig:15} shows the distribution of fair values of the instrument, plus the coupons in the respective path obtained so far after five years and ten years in the future. The histogram gives an idea about the frequency distribution and shape of a set of 10000 values for the PRIIP.

\section{Conclusion}
This paper presents a detailed error analysis of the model order reduction framework developed in \cite{binder2020} for the financial risk analysis. Three major sources of errors considered are the discretization error, the sampling error, and the model order reduction error. We have developed error estimations and their quantification methods for each error. Based on these error estimators, we have designed algorithms to select optimal grid size, suitable training parameters, and the appropriate reduced dimension. The algorithms are designed such that the model order reduction framework can be used as a black box for which the user only requires to provide inputs and get the desired output.

Another major contribution of this paper is the implementation of a sensitivity analysis to rank the contribution of each numerical error. We have developed a methodology to address both deterministic (e.g., discretization error) and stochastic errors (sampling error). The Sobol variance-based sensitivity analysis is used, which computes local as well as global indices. The global sensitivity analysis revealed that the sampling error is most sensitive to the final output. Such information advised us to allocate more computational resources to the sampling algorithm. This ultimately generated a suitable reduced basis as demonstrated by the obtained results.

We have tested the designed algorithms for a numerical example of a puttable steepener under the two-factor Hull-White model. The reduced model is obtained such that the total error $\epsilon_T = 6.479\times 10^{-4}$ is less than the user-defined tolerance $e_{tol} = 10^{-3}$. 
As $\epsilon_T < e_{tol}$, we can conclude that the full model valuations for the adaptively selected $\ell = 10$ parameter groups are enough, and the rest of the parameter groups can be solved inexpensively using the reduced model of dimension $d=10$. The reduced model approach to solve the entire parameter space was at least $8-10$ times faster than the full model approach. We also noticed that the randomized singular value decomposition provides an excellent speedup compared to the basic singular value decomposition.

We conclude that the new model order reduction approach shows many potential applications in the historical or Monte Carlo value at risk calculations.
\section*{Acknowledgement}
This project has received funding from the European Union’s Horizon 2020 research and innovation program under the Marie Skłodowska-Curie Grant Agreement No. 765374.\\
We thank our colleagues Michael Schwaiger and Diana Hufnagl from MathConsult who provided insight and expertise that greatly assisted the research.
\bibliography{references}
\end{document}